\documentclass[12pt]{amsart}
\usepackage{mathrsfs}
\usepackage{amscd}
\usepackage{amsmath}
\usepackage{latexsym}
\usepackage{amsfonts}
\usepackage{amssymb}
\usepackage{amsthm}
\usepackage{graphicx}

\def\R{\mathbb{R}}

\def\T{\mathbb{T}}
\def\C{\mathbb{C}}

\newtheorem{thm}{Theorem}[section]
\newtheorem{lem}{Lemma}[section]
\newtheorem{prop}{Proposition}[section]

\newtheorem{remark}{Remark}[section]
\newtheorem{definition}{Definition}[section]

\makeatletter
\newcommand{\Extend}[5]{\ext@arrow0099{\arrowfill@#1#2#3}{#4}{#5}}
\makeatother

\begin{document}

\setcounter{page}{1}

\title[Global rough solution for critical gKdV]{The low regularity global solutions for the critical generalized KdV equation}
\author{Changxing Miao}
\address{Institute of Applied Physics and
Computational Mathematics,  P. O. Box 8009,\ Beijing,\ China,\
100088,}
\email{miao\_changxing@iapcm.ac.cn}

\author{Shuanglin Shao}
\address{Institute for Mathematics and its Applications, University of Minnesota, Minneapolis, MN 55455}
\email{slshao@ima.umn.edu}

\author{Yifei Wu}
\address{Department of Mathematics, South China University of Technology,
Guangzhou, Guangdong, China, 510640,}
\email{yerfmath@yahoo.cn}

\author{Guixiang Xu}
\address{Institute of
Applied Physics and Computational Mathematics,  P. O. Box 8009,\
Beijing,\ China,\ 100088, }
\email{xu\_guixiang@iapcm.ac.cn}


\date{\today}


\maketitle

\begin{abstract}\noindent
We prove that the Cauchy problem of the mass-critical generalized KdV equation is globally well-posed in Sobolev spaces $H^s(\R)$ for $s>6/13$. Of course, we require that the mass is strictly less than that of the ground state in the focusing case. The main approach is the ``I-method" together with the multilinear correction analysis. Moreover, we
use some ``partially refined" argument to lower the upper control of the multiplier in the resonant interactions. The result improves the previous works of Fonseca, Linares, Ponce (2003) and Farah (2009).
\end{abstract}

\section{Introduction}
In this paper, we consider the global well-posedness of the Cauchy problem for the mass-critical generalized Korteweg-de Vries equation (gKdV):
\begin{eqnarray}\label{gkdv}
  && \partial_{t}u+\partial_{x}^{3}u=\mu\partial_x(u^5),
     \qquad u:\R\times [0,T]\mapsto \R,\\
  && u(x,0)=u_0(x)\in H^s(\R),\label{1.2}
\end{eqnarray}
where $\mu=\pm 1$, $H^s(\R)$ denotes the usual inhomogeneous Sobolev space of order $s$. When $\mu =1$, the equation \eqref{gkdv} is called `` defocusing", while when $\mu=-1$ it is called ``focusing". The equation \eqref{gkdv} is mass-critical since the scaling $$u(x,t)\to \lambda^{-1/2}u(x/\lambda,t/\lambda^3), \,\lambda>0, $$
leaves both the equation and the mass $\int_{\mathbb{R}} |u(x,t)|^2 dx $ invariant. It is well-known that
\eqref{gkdv} belongs to a family of equations,
$$
\partial_t u +\partial_x^3 u  =\mu \partial_xu^{p}, \quad u:
\mathbb{R}\times [0, T]\longmapsto  \mathbb{R},\\
$$
where $p\ge 3$.

The Cauchy problem \eqref{gkdv}-\eqref{1.2} was shown by Kenig, Ponce, Vega \cite{KPV3} to be locally well-posed in $H^s(\R)$ for $s\geq 0$, see also \cite{K} for $s>3/2$. In other words, for any $u_0\in H^s(\R)$, there exists a positive time $T=T(\|u_0\|_{H^s})$ (when $s=0$, $T$ also depends on the profile of the initial data $u_0$), such that the solution to \eqref{gkdv}-\eqref{1.2} exists and is unique in a certain Banach space of functional $X\subset C\bigl([0,T];H^s(\R)\bigr)$; moreover, the solution map is continuous from $H^s(\R)$ to $C([0,T];H^s(\R))$. It appears that the index $s=0$ is sharp since there are examples to show that the critical gKdV equation is ill-posed for $s<0$, see \cite{BKPSV}. If the lifetime of the solution $T$ can be taken arbitrarily large, we say that \eqref{gkdv}-\eqref{1.2} is globally well-posed.

As is well-known, the local solution to equation \eqref{gkdv}-\eqref{1.2} enjoys the \emph{mass conservation law},
\begin{equation}\label{MCL}
M(u(t))\equiv\displaystyle\int_{\R} |u(x,t)|^2 \,dx =M(u_0),
\end{equation}
and the $H^1-$solution enjoys the \emph{energy conservation law},
\begin{equation}\label{ECL}
E(u(t))\equiv\displaystyle\int_{\R} \frac{1}{2}|\partial_x
u(x,t)|^2+\frac{\mu}{6}|u(x,t)|^6 \,dx =E(u_0).
\end{equation}
Hence, an immediate conclusion on global wellposedness for $H^1$-initial data follows from the local theory in \cite{KPV3} and the equation \eqref{ECL} above in the defocusing case. In the focusing case, the same conclusion holds
under the condition $\|u_0\|_{L^2}< \|Q\|_{L^2}$ by the sharp
Gagliardo-Nirenberg inequality (see \cite{W}),
\begin{equation}\label{GN}
\|u\|_{L^6}^6\leq 3(\|u\|_{L^2}/\|Q\|_{L^2})^4\|\partial_xu\|_{L^2}^2,
\end{equation}
where $Q=[3\mbox{sech}^2(2x)]^{\frac{1}{4}}$ is the \emph{ground state solution} to
the elliptic equation
$$
\partial_{xx}Q+Q^5=Q.
$$
Moreover, the local theory in \cite{KPV3} implies the global well-posedness in $L^2$ when the initial data has sufficiently small $L^2$ norm.  However, unlike in the case of $H^1(\R)$, where the equation
\eqref{gkdv} is ``subcritical" with respect to the regularity of initial data, the usual iteration argument involving the Strichartz estimates and the \emph{mass conservation law} will not yield the global wellposedness directly for large $L^2 $ data. So the question of $L^2$-global wellposedness and scattering is regarded as an open conjecture in
the field; it is far from resolution, despite much recent progress \cite{KKSV, Shao-08}.

Therefore, a natural question arises: what is the least $s_0>0$ such that for $s>s_0$, if $u_0\in H^s(\R)$, the solution to \eqref{gkdv}-\eqref{1.2} is globally well-posed? The question is partly plausible in light of the recent exciting progresses in nonlinear dispersive equations such as nonlinear Schr\"odinger equations
(NLS), nonlinear wave equations (NLW), etc. They are made possibly by the well-known strategies: Bourgain's ``Fourier truncation method" in \cite{Bourgain} and the ``I-method" by I-team (Colliander, Keel, Staffilani, Takaoka, Tao) partially inspired by the former, see e.g.,  \cite{CKSTT-01-KDV}, \cite{CKSTT-03-KDV}. The ``Fourier truncation method" works well provided that there is some smoothing effect arising from the non-linearity, while ``I-method" can still work in the case that there are derivatives in the non-linearity and such smoothing is not available, and often the latter gives a sharper result (see \cite{KT98, KT2} for a discussion). We do not intend to survey these two methods, but we refer readers to \cite{CKSTT-01-KDV} for a nice discussion on Bourgain's high/low trick in ``Fourier truncation method" and in Section 2 and 3 of that paper readers can also find an introduction and an example of applications of the first-generation ``I-method" and second generation ``I-method" to KdV and mKdV equations. For a textbook treatment, we refer readers to \cite[Chapter 3.9]{Tao:2006-CBMS-book}. For the recent developments of ``I-method", we can refer to \cite{CGT07, CKSTT-02-NLS, CKSTT-04, CKSTT-08, CR08, DPST07, DPST08,
FG07, VZ09} on the applications in the context of nonlinear Schr\"{o}dinger equation (NLS), refer to \cite{CKSTT-01-DNLS, CKSTT-02-DNLS} on the applications in the context of Schr\"{o}dinger equation with derivative (DNLS), refer to \cite{CKSTT-01-KDV, CKSTT-03-KDV, CST99, GPS07, Zhang} on the applications in the context of gKdV equations.

The global well-posedness of \eqref{gkdv}-\eqref{1.2} below the energy space $H^1$ was considered by Farah \cite{F}, Fonseca, Linares and Ponce \cite{FLP}. The authors in \cite{FLP} proved the global existence in $H^s(\R)$ for $s>3/4$ by appllying Bourgain's ``Fourier truncation method". This was improved very recently in \cite{F}, which lower the index to $s>3/5$ by I-method introduced in \cite{CKSTT-02-NLS}. The condition $\|u_0\|< \|Q\|_{L^2}$ is imposed for both results in the focusing case.

Our main result in this paper is the following improvement.
\begin{thm}\label{thm:main}
The Cauchy problem \eqref{gkdv}-\eqref{1.2} is globally well-posed in $H^s(\R)$ for $s>6/13$ when $\mu=1$. The same conclusion holds under the assumption $\|u_0\|< \|Q\|_{L^2}$ when $\mu=-1$.
\end{thm}

To prove this theorem, we will follow the general scheme of ``I-method" by adding a ``correction-term" to the first modified energy $E(Iu)$ as in \cite{CKSTT-02-DNLS, CKSTT-03-KDV} and using the multilinear correction analysis. If we added it in a naive way, the multiplier introduced in order to obtain the second modified energy is singular in the sense that its $L^\infty$ norm is infinity for a set of nonzero measure, see Section \ref{sec:I-method} of this paper for an exact description. This difficulty is also noted in \cite[Proposition 3.1]{F}. To get out of it, our approach here follows along similar lines of refining the modified energy by performing a resonant decomposition to the singular multiplier as in \cite{CKSTT-08} by Colliander, Keel, Staffilani, Takaoka and Tao,
see also \cite{BDGS, Bourgain2}.  More precisely we will split the multiplier $M_6$ arising from the derivative of  the first modified energy into two parts:
 $$M_6:=\bar{M}_6+\tilde{M}_6$$
under the principle that $\bar{M}_6$ contains some low-frequency terms, which is referred to as ``resonant term",  and $\tilde{M}_6$ contains the rest, which is referred to as the ``non-resonant term". For $\tilde{M}_6$, we will use a point-wise estimate (in $t$) and reduce it to an error term in the final bootstrap argument, see
Lemma \ref{le:pointwise-bound} and the argument in Section \ref{sec:final-argument}. For $\bar{M}_6$, we perform a careful multilinear analysis by using $X_{s,b}$-type estimates.

Our key point for such an improvement $s> {6}/{13}$ is due to a better energy increment bound,  $N^{-\frac 72+}$, see the statements in Theorem \ref{thm:main-2}. This improves the previous estimate, $N^{-2+}$ in \cite{F}, which gives $s>{3}/{5}$.

We will focus on the focusing equation \eqref{gkdv} under the assumption that $\|u_0\|_2<\|Q\|_2$. This assumption guarantees that the kinetic energy in \eqref{ECL} is comparable to the energy thanks to the sharp Gagliardo-Nirenberg inequality \eqref{GN}. The same analysis will go through the defocusing case in the same way.

The paper is organized as follows. In Section \ref{sec:notation}, we introduce some notations and state some preliminary estimates that will be used throughout this paper. In Section \ref{sec:I-method}, we establish a variant of local well-posedness theory, set up the I-method. In Section \ref{sec:main-arguments}, we establish some fixed time bound for the error term and obtain an upper bound on the increment of the new modified energy. In Section \ref{sec:final-argument}, we prove the global well-posedness in Theorem \ref{thm:main}.

\section{Notations and some preliminary estimates}\label{sec:notation}
\subsection{Notations.} We use $A\lesssim B$ or $B\gtrsim A$ to denote the statement that $A\leq CB$ for some positive constant $0<C<\infty$ which does not depend on the functions but may vary from line to line. We use the notation $A\sim B$ whenever $A \lesssim B$ and $B\lesssim A$. If the constants appearing in $\lesssim$ or $\gtrsim$ depend upon some additional parameters, we will indicate them with subscripts; for example, $A \lesssim_\epsilon B$ denotes the assertion that $A \le C_\epsilon B$  for some positive constant $C_\epsilon$ depending on $\epsilon$; similarly for $A \sim_\epsilon B$, etc.

We use $A\ll B$, or sometimes $A=o(B)$ to state the statement $A\leq C^{-1}B$ for a sufficiently large constant $C>0$.
The notation $a+$ denotes $a+\epsilon$, and $a-$ for $a-\epsilon$ for arbitrarily small exponents $\epsilon >0$, and allow the implied constants in $\lesssim $ notation to depend on $\epsilon$.

We also set
$$\langle\cdot\rangle=(1+|\cdot|^2)^{1/2}, \, D_x^\alpha=(-\partial^2_x)^{\alpha/2}, \,J_x^\alpha=(1-\partial^2_x)^{\alpha/2}.$$
We use $\|f\|_{L^p_xL^q_t}$ to denote the mixed norm $\Big(\displaystyle\int\|f(x,\cdot)\|_{L_t^q}^p\
dx\Big)^{\frac{1}{p}}$.

Now we record some definitions. For $s,\,b\in \R $, we define the Bourgain space $X_{s,\,b}$ to be the closure of the Schwartz class under the norm
\begin{equation}\label{X}
\|u\|_{X_{s,b}}\equiv \left(\int\!\!\!\!\int
\langle\xi\rangle^{2s}\langle\tau-\xi^3\rangle^{2b}|\tilde{u}(\xi,\tau)|^2\,d\xi
d\tau\right)^{\frac{1}{2}},
\end{equation}where $\tilde{u}$ denotes the space-time Fourier transform of $u$ defined by
$$\tilde{u}(\xi,\tau)=\iint e^{-i(x\xi+t\tau)} u(x,t) dxdt;$$
similarly we denote by $\widehat{f}$ the Fourier transform of $f(x,t)$ in the spatial variable. For any interval $\Omega$, we define $X_{s,b}^{\Omega}$ to be the restriction of $X_{s,b}$ on $\R\times\Omega$ with the norm
\begin{equation}\label{X1}
\|u\|_{X_{s,b}^\Omega}=\inf\{\|U\|_{X_{s,b}}:U|_{t\in\Omega}=u|_{t\in\Omega}\}.
\end{equation}
When $\Omega=[-\delta,\delta]$, we will write $X_{s,b}^\Omega$ as $X_{s,b}^\delta$.

Let $0<s<1$ and $N\gg 1$ be fixed. The Fourier multiplier operator $I_{N,s}$ is defined by
\begin{equation}\label{I}
\widehat{I_{N,s}u}(\xi)=m_{N,s}(\xi)\widehat{u}(\xi),
\end{equation}
where the multiplier $m_{N,s}(\xi)$ is a smooth, monotone and radial function satisfying $0<m_{N,s}(\xi)\leq 1$ and
\begin{equation}\label{m}
m_{N,s}(\xi)=\biggl\{
\begin{array}{ll} 1,&|\xi|\leq N,\\ \bigl(\frac {N}{|\xi|}\bigr)^{1-s},&|\xi|>2N.
\end{array}
\end{equation}
Sometimes we denote $I_{N,s}$ and $m_{N,s}$ by $I$ and $m$ respectively if there is no confusion.

\begin{remark}
The operator $I_{N,s}$ maps $H^s(\R)$ into $H^1(\R)$ with equivalent norms for any $s<1$. More precisely, there exists some positive constant $C$ such that
\begin{equation}\label{II}
C^{-1}\|u\|_{H^s}\leq \|I_{N,s}u\|_{H^1} \leq CN^{1-s}\|u\|_{H^s}.
\end{equation} Moreover, $I_{N,s}$ can be extended to a map (still denoted by $I_{N,s}$) from $X_{s,b}$ to $X_{1,b}$ which satisfies that for any $s<1,b \in \R$,
$$ C^{-1}\|u\|_{X_{s,b}}\leq\|I_{N,s}u\|_{X_{1,b}} \leq CN^{1-s}\|u\|_{X_{s,b}}. $$
\end{remark}

\subsection{Preliminary estimates} We state some preliminary estimates which will be used throughout the paper. We start with some well-known Strichartz estimates, see e.g., \cite{KPV4, KPV3}.
\begin{lem}\label{le:strichartz-esimate-1}
For $u\in X_{0,\frac{1}{2}+}$, we have
\begin{equation}
\|D_xu\|_{L^\infty_x L^2_t}+\left\|D_x^{-\frac{1}{4}}u\right\|_{L^4_x L^\infty_t}+
\|u\|_{L^5_x L^{10}_t}\lesssim \|u\|_{X_{0,\frac{1}{2}+}}.\label{XE1}
\end{equation}
\end{lem}
By Sobolev's embedding, we have the following estimate.
\begin{lem}\label{le:strichartz-esimate-2}
 For $u\in X_{\frac{1}{2}+,\frac{1}{2}+}$, we have
\begin{equation}
\|u\|_{L^\infty_{x,t}} \lesssim \|u\|_{X_{\frac{1}{2}+,\frac{1}{2}+}}.\label{XE2}
\end{equation}
\end{lem}
We recall the following bilinear estimate from \cite[Corollary 3.2]{G} and \cite[Lemma 2.1]{F}.
\begin{lem}\label{le:strichartz-esimate-3}
For any $f_1,f_2\in X_{0,\frac{1}{2}+}$ supported on the frequencies $\{|\xi_i|\sim N_i\}$, $i=1,2$. If $|\xi_1|\sim |\xi_1-\xi_2|\sim |\xi_1+\xi_2|$ for all $\xi_i\in \operatorname{Supp}\widehat{f_i}$, $i=1,2$, then
\begin{equation}
\left\|(D_xf_1)f_2\right\|_{L^2_{x, t}} \lesssim
\|f_1\|_{X_{0,\frac{1}{2}+}}\,\|f_2\|_{X_{0,\frac{1}{2}+}}. \label{IE1}
\end{equation}
\end{lem}
For the sake of completeness, we provide a proof by using Plancherel's theorem. 
\begin{proof} To prove \eqref{IE1}, it suffices to prove
\begin{equation}\label{eq-53}
\|D_xe^{-t\partial_x^3}\phi_1 e^{-t\partial_x^3}\phi_2\|_{L^2_{t,x}} \lesssim \|\phi_1\|_{L^2}\|\phi_2\|_{L^2},
\end{equation}
where $|\xi_1|\sim |\xi_1-\xi_2|\sim |\xi_1+\xi_2|$ for all $\xi_i\in \operatorname{Supp}\widehat{\phi_i}$. We write
$$D_xe^{-t\partial_x^3}\phi_1 e^{-t\partial_x^3}\phi_2 =\iint e^{ix(\xi_1+\xi_2)+it(\xi_1^3+\xi_2^3)}|\xi_1|\widehat{\phi_1}\widehat{\phi_2} d\xi_1d\xi_2.$$
We change variables as follows, $a:=\xi_1+\xi_2$ and $b:=\xi_1^3+\xi_2^3$; then the Jocabian 
$$J:=\left|\frac {\partial(a,b)}{\partial(\xi_1,\xi_2)}\right|\sim |\xi_1^2-\xi_2^2|=|(\xi_1+\xi_2)(\xi_1-\xi_2)|\sim|\xi_1|^2$$
by the assumption on the Fourier supports of $\phi_i$, $i=1,2$. Then we apply Plancherel's theorem to the left hand side of \eqref{eq-53} followed by a changing of variables back, we see that it is bounded by
\begin{equation}\label{eq-54}
\lesssim \left(\iint \frac {\left(|\xi_1|\widehat{\phi_1}\widehat{\phi_2}\right)^2}{J}d\xi_1 d\xi_2\right)^{1/2} \lesssim \|\phi_1\|_{L^2}\|\phi_2\|_{L^2}.
\end{equation} 
This proves \eqref{eq-53}, and hence Lemma \ref{le:strichartz-esimate-3}.
\end{proof}

\section{I-method and the Multilinear Estimates}\label{sec:I-method}

\subsection{A Variant Local Well-posedness}
In this subsection, we will establish a variant local well-posedness result.
\begin{prop}\label{prop:modified-local}
Let $s>0$, then Cauchy problem \eqref{gkdv}-\eqref{1.2} is locally well-posed for the initial data $u_0$ with $I_{N,s}u_0\in H^1(\R)$. Moreover, the solution $u$ exists on the interval $[0,\delta]$ with the lifetime
\begin{equation}
 \delta\sim\|I_{N,s}u_0\|^{-\mu}_{H^1}\label{delta}
\end{equation}
for some $\mu>0$, and
\begin{equation}
\|I_{N,s}u(t)\|_{X_{1,\frac{1}{2}+}^\delta}\lesssim \|I_{N,s}u_0\|_{H^1}.\label{LSE}
\end{equation}
\end{prop}
It can be established by a standard iteration argument; we present it here for sake of completeness;  see also \cite[Theorem 5.1]{F} for $s>\frac12$.

 \begin{proof} The proof proceeds by the usual fixed point argument on the space $X_{1,\frac12+}^I(J)$.
 By Duhamel's principle, Lemma 2.11 and Proposition 2.12 in \cite{Tao:2006-CBMS-book}, we have
 \begin{equation*}
 \aligned \big\|Iu\big\|_{X_{1,\frac12+}(J)} &= \Big\|S(t)\big(Iu_0\big)
 +\int^t_0 S(t-s)I\big( u^4\partial_x u \big)(s)ds \Big\|_{X_{1,\frac12+}(J)}\\
 & \lesssim \big\|Iu_0\big\|_{H^1}
 + \big\| I\big( u^4\partial_x u \big) \big\|_{X_{1, -\frac12+}(J)}\\
 & \lesssim \big\|Iu_0\big\|_{H^1}+ \delta^{\epsilon} \big\| I\big(
 u^4\partial_x u \big) \big\|_{X_{1, -\frac12++}(J)} \endaligned
 \end{equation*}
 for sufficiently small $\epsilon>0$, where $a++=a+2\epsilon$. Hence it suffices to show that
 \begin{equation*}
 \aligned \big\| I\big( u^4\partial_x u \big) \big\|_{X_{1, -\frac12++}} &\lesssim \big\| Iu  \big\|^5_{X_{1, \frac12+}}.
 \endaligned
 \end{equation*}
 Using Lemma 12.1 in \cite{I-team:2004:multilinear-estimates-for-periodic-KdV-and-applications} or the argument of Lemma 5.2 in \cite{CKSTT-01-DNLS}, we only need to prove that
 \begin{equation}\label{nonlinear}
 \aligned \big\|  u^4\partial_x u  \big\|_{X_{s, -\frac12++}} & \lesssim \big\| u
 \big\|^5_{X_{s, \frac12+}} \quad \text{for}\quad 0<s<1. \endaligned
 \end{equation}
 Indeed, from Lemma \ref{le:strichartz-esimate-1}, we have
 \begin{equation*}
 \aligned \big\|f\big\|_{L^5_xL^{10}_t} \lesssim \big\|f\big\|_{X_{0,\frac{1}{2}+}}. \endaligned
 \end{equation*}
 which, interpolating with $\big\|f\big\|_{L^{\infty}_xL^{\infty}_t}
 \lesssim \big\|f\big\|_{X_{\frac12+,\frac12+}} $, implies that
 \begin{equation*}
 \aligned \big\|f\big\|_{L^{5+}_xL^{10+}_t} \lesssim \big\|f\big\|_{X_{0+,\frac{1}{2}+}}. \endaligned
 \end{equation*}
 In addition, by duality of $\big\|f\big\|_{L^5_xL^{10}_t} \lesssim \big\|f\big\|_{X_{0,\frac{1}{2}+}}$, we also have
 \begin{equation*}
 \aligned \big\|f\big\|_{X_{0,-\frac{1}{2}-}} \lesssim \big\|f\big\|_{L^{\frac{5}{4}}_xL^{\frac{10}{9}}_t}, \endaligned
 \end{equation*}
 which, interpolating with $\big\|f\big\|_{X_{0,0}} = \big\|f\big\|_{L^{2}_xL^{2}_t}$, implies that
 \begin{equation*}
 \aligned \big\|f\big\|_{X_{0,-\frac{1}{2}++}} \lesssim
 \big\|f\big\|_{L^{\frac{5}{4}+}_xL^{\frac{10}{9}+}_t}. \endaligned
\end{equation*}
 Hence, by the fractional Leibniz rule (Principle A.5 in \cite{Tao:2006-CBMS-book}) and Lemma \ref{le:strichartz-esimate-1}, we have
 \begin{equation*}
 \aligned \big\|  u^4\partial_x u  \big\|_{X_{s, -\frac12++}} & \lesssim
 \big\|u^3\langle \nabla \rangle^{s} u \partial_x u \big\|_{L^{\frac{5}{4}+}_xL^{\frac{10}{9}+}_t} + \big\|u^4\langle \nabla \rangle^{s} \partial_x u \big\|_{L^{\frac{5}{4}+}_xL^{\frac{10}{9}+}_t}\\
 & \lesssim \big\| u \big\|^3_{L^{5+}_xL^{10+}_t} \big\|\langle
 \nabla \rangle^{s} u \big\|_{L^5_xL^{10}_t} \big\| \partial_x u
 \big\|_{L^{\infty}_xL^{2}_t} +  \big\| u \big\|^4_{L^{5+}_xL^{10+}_t} \big\|\langle \nabla \rangle^{s} \partial_x u \big\|_{L^{\infty}_xL^{2}_t} \\
 &\lesssim \big\| u \big\|^3_{X_{0+, \frac12+}} \big\|\langle \nabla
 \rangle^{s} u \big\|_{X_{0,\frac12+}} \big\| u \big\|_{X_{0,\frac12+}}
 + \big\| u \big\|^4_{X_{0+,\frac12+}} \big\|\langle \nabla \rangle^{s} u \big\|_{X_{0,\frac12+}}\\
 & \lesssim \big\| u \big\|^5_{X_{s, \frac12+}}\endaligned
 \end{equation*}
 for $s>0$. This completes the proof of Proposition \ref{prop:modified-local}.
 \end{proof}

\subsection{I-method and modified energy}
From now on, we take $\mu=-1$ and $u$ be the real-valued solution of \eqref{gkdv}-\eqref{1.2} throughout the paper.

First we record the classical set-up for the ``I-method", see Section 2 in \cite{CKSTT-03-KDV} or Section 3 in \cite{CKSTT-08}. Given a smooth tempered symbol $M_k(\xi_1,\cdots,\xi_k)$ defined on the hyperplane with the the push-forward Lebesgue measure $d\xi_1\cdots d\xi_{k-1}$,
\begin{equation}\label{Gamma_n}
\Gamma_k:=\left\{(\xi_1,\cdots, \xi_k):\xi_1+\cdots+\xi_k=0\right\},
\end{equation}
we define the quantity
\begin{equation*}
  \Lambda_k(M_k):=\displaystyle\int_{\Gamma_k}M_k(\xi_1,\cdots,\xi_k)\prod_{j=1}^k\widehat{u}(\xi_j,t)\,d\xi_1\cdots d\xi_{k-1}.
\end{equation*}
Then by using the equation \eqref{gkdv} and a direct computation, we have the following differentiation formula.
\begin{lem}[Differentiation formula]\label{le:dLMD}Let $\Lambda_k$ and $M_k$ be defined as above. Then
\begin{equation}\label{dLMD}
\dfrac{d}{dt}\Lambda_k(M_k)=\Lambda_k(M_k\alpha_k)-ik\Lambda_{k+4}(M_k(\xi_1,\cdots,\xi_{k-1},\xi_k+\cdots+\xi_{k+4})
(\xi_k+\cdots+\xi_{k+4})),
\end{equation}
where $\alpha_k$ is the symbol defined by $$\alpha_k:=i(\xi_1^3+\cdots+\xi_k^3).$$
\end{lem}

We define the ``first-generation" modified energy by
\begin{equation}
E^1_I(u(t)):=\dfrac{1}{2}\|\partial_xIu(t)\|_{L^2}^2-\dfrac{1}{6}\|Iu(t)\|_{L^6}^6.
\end{equation}
Then it follows from the Fourier inversion formula that
\begin{equation}\label{E1}
E^1_I(u(t))=\Lambda_2(\sigma_2)+\Lambda_6(\sigma_6)
\end{equation} where $\sigma_2$ and $\sigma_6$ are symbols defined by
$$
\sigma_2:=-\dfrac{1}{2}m(\xi_1)m(\xi_2)\xi_1\xi_2;\quad
\sigma_6:=-\dfrac{1}{6}m(\xi_1)\cdots m(\xi_6).
$$
Furthermore by the differentiation formula (\ref{dLMD}), we have
\begin{equation}\label{dE1}
\begin{split}
\dfrac{d}{dt}E^1_I(u(t))=&\; \Lambda_2(\sigma_2\alpha_2)
-2i\Lambda_{6}(\sigma_2(\xi_1,\xi_{2}+\cdots+\xi_{6})(\xi_{2}+\cdots+\xi_{6}))+\Lambda_6(\sigma_6\alpha_6)\\
&\;-6i\Lambda_{10}\Big(\sigma_6(\xi_1,\cdots,\xi_{5},\xi_{6}+\cdots+\xi_{10})
(\xi_{6}+\cdots+\xi_{10})\Big)\\
=&\Lambda_{6}\bigl(-2i\sigma_2(\xi_1,\xi_{2}+\cdots+\xi_{6})(\xi_{2}+\cdots+\xi_{6})+\sigma_6\alpha_6\bigr)\\
&\;+\Lambda_{10}\Big(-6i\sigma_6(\xi_1,\cdots,\xi_{5},\xi_{6}+\cdots+\xi_{10})
(\xi_{6}+\cdots+\xi_{10})\Big)\\
=:&\; \Lambda_{6}(M_6)+\Lambda_{10}(M_{10}),
\end{split}
\end{equation}
Note that the first term vanishes because $\alpha_2$ vanishes on $\Gamma_2$. Here by a similar symmetrization consideration as in \cite[Section 2 and 3]{CKSTT-03-KDV} or \cite[Section 3]{CKSTT-08}, we have
\begin{equation}\label{eq:sym-reduction}
 \Lambda_6(M_6, u)=\Lambda_6([M_6]_{sym}), \, \Lambda_{10}(M_{10}, u)=\Lambda_{10}([M_{10}]_{sym});
 \end{equation}
Since $\alpha_6$ and $\sigma_6$ are already symmetric with respect to the group $S_k$, the group of all permutations on $k$ objects. we have
\begin{equation}\label{eq:sym-formula}
\begin{split}
[M_6]_{sym}&=-2i[\sigma_2(\xi_1,\xi_{2}+\cdots+\xi_{6})(\xi_{2}+\cdots+\xi_{6})]_{sym}
 +\sigma_6\alpha_6\\
 &=:M_6^1+M_6^2;\\
[M_{10}]_{sym}&=-6i[\sigma_6(\xi_1,\cdots,\xi_{5},\xi_{6}+\cdots+\xi_{10})
(\xi_{6}+\cdots+\xi_{10})]_{sym}.
\end{split}
\end{equation}
For readers' convenience, we record the definition of the symmetrization of a multiplier from \cite[Definition 1]{CKSTT-03-KDV}.
\begin{definition} A $k$-multiplier is a function $m: \R^k\to \C$. A $k$-multiplier is symmetric if $m(\xi)= m\bigl(g(\xi)\bigr)$ for all $g\in S_k$. The symmetrization of a $k$-multiplier $m$ is the multiplier
$$[m]_{sum}(\xi):=\frac {1}{k!}\sum_{g\in S_k} m\bigl(g(\xi)\bigr).$$
\end{definition}
\begin{remark}\label{re:sym-consideraton}
An example of a symmetric 2-multiplier is $\sigma_2=-\frac {1}{2}m(\xi_1)m(\xi_2)\xi_1\xi_2$ defined above, and by an explicit computation
\begin{align*}
M_6^1&=-2i[\sigma_2(\xi_1,\xi_{2}+\cdots+\xi_{6})\bigl(\xi_{2}+\cdots+\xi_{6})]_{sym}\\
&=\frac {i}{6}(m^2(\xi_1)\xi_1^3+\cdots +m^2(\xi_6)\xi_6^3\bigr).\\
[M_{10}]_{sym}&=-6i[\sigma_6(\xi_1,\cdots,\xi_{5},\xi_{6}+\cdots+\xi_{10})
(\xi_{6}+\cdots+\xi_{10})]_{sym}\\
&=C\sum_{\{a,\cdots,j\}= \{1,\cdots, 10\}} m(\xi_a) m(\xi_b)
m(\xi_c) m(\xi_d) m(\xi_e) m(\xi_f+\cdots +\xi_j)\\
&\qquad\qquad \times
(\xi_f+\cdots +\xi_j)
\end{align*}for some explicit nonzero constant $C$.
\end{remark}

\begin{remark}[Two convenient reductions]\label{re:two-reductions}
There are two well-known reductions which we will use throughout the rest of the paper.
\begin{enumerate}
\item [(1)] By symmetrization, the first reduction is that we could order the magnitudes of $\xi_i$: for example, assume that $|\xi_1|\ge\cdots \ge |\xi_6|$ in estimating $\Lambda_6$.
\item [(2)] In various 6-linear or 10 linear estimates below, we may make a Littlewood-Paley decomposition and restrict attention to the contribution arising in $|\xi_i|\sim N_i$ where $N_i$ is dyadic; from the discussion above, we may assume that $N_1\ge \ldots \ge N_{10}$. If $N_1\ll N$,
    $$M_6=0, \,\text{ and } M_{10}=0 \Rightarrow E^1_I(u(t)) \text{ is conserved for all time.}$$
    Then global wellposedness for this solution $u$ would follow. Hence without loss of generality, we will take $N_1\gtrsim N$, then the support information of $\Gamma_6$ or $\Gamma_{10}$ will give $|N_2|\gtrsim N$, and hence $N_1\sim N_2\gtrsim N$, which is our second reduction.
    \end{enumerate}
\end{remark}

Now we elaborate the difficulty if following  a direct analogous reasoning as in \cite{CKSTT-03-KDV}, which forces us to think of an alternative by introducing a resonant decomposition on the multiplier: From the expression for $\frac {d}{dt}E^1_I(u(t))$ in \eqref{dE1}, a natural choice of the second modified energy would be
$$E^2_I(u(t))=\Lambda_6(\check{\sigma}_6)+E^1_I(u(t))$$
with the choice
$$\check{\sigma}_6:=-M_6/\alpha_6.$$ But unfortunately, $\check{\sigma}_6$ is singular and hence is unfavorable. This is in contrast with the cases for KdV or mKdV \cite{CKSTT-03-KDV}, where one can take advantage of the complete integrability for these equations.

As forecasted in the Introduction, to overcome this difficulty, the strategy here is to split $M_6$
into two parts: $$M_6=\bar{M}_6+\tilde{M}_6,$$
where $\bar{M}_6$ and $\tilde{M}_6$ are defined to be ``resonant" and ``non-resonant" parts, respectively. Let us motivate the choice of the non-resonant set. Roughly speaking, what we expect is that,
\begin{itemize}
\item \textit{either} $M_6$ \textit{is controlled by a low frequency term,}
\item \textit{ or non-resonant occurs, i.e., } $|M_6|\lesssim |\alpha_6|$.
\end{itemize}

Suppose $|\xi_1|\geq |\xi_2|\gtrsim N \gg |\xi_3|, \ldots, |\xi_6|$.  Then
in the resonant case, one may find $\alpha_6=0$ and thus
$$
\xi_1^3+\xi_2^3=-\bigl(\xi_3^3+\xi_4^3+\xi_5^3+\xi_6^3\bigr).
$$
Then we have a coarse estimate
\begin{equation}\label{eq-48}
|\xi_1^3+\xi_2^3|\lesssim|\xi_3^3+\xi_4^3+\xi_5^3+\xi_6^3|.
\end{equation}
Together with the information $|\xi_2|\gtrsim N\gg |\xi_3|$, \eqref{eq-48} implies:
\begin{equation}\label{eq-49}
|M_6|\lesssim|\xi_3^3+\xi_4^3+\xi_5^3+\xi_6^3|,
\end{equation} as $\big|m^2(\xi_1)\xi_1^3+m^2(\xi_2)\xi_2^3 \big| \lesssim |\xi_1^3+\xi_2^3|$  by the mean value theorem in Lemma \ref{le:mvt}. This shows that $M_6$ is bounded by a lower frequency term as expected. So if we take the contrapositive to \eqref{eq-48}, i.e.,
$$
\big|\xi_1^3+\xi_2^3\big|\gg\big|\xi_3^3+\xi_4^3+\xi_5^3+\xi_6^3\big|,
$$
we are in the non-resonant case. If Suppose that $|\xi_1|\geq|\xi_2|\geq|\xi_3|\gtrsim N \gg
|\xi_4|\geq|\xi_5|\geq|\xi_6|$, one may find it is always non-resonant, this motivates the choice of $\Omega_2$. Lastly if we are in the case where $|\xi_1|\geq\cdots \geq|\xi_4|\gtrsim N \gg|\xi_5|\geq|\xi_6|$, then the motivation of the choice of $\Omega_3$ is similar to $\Omega_1$ but a little more complicated.

Let us define ``non-resonant" sets. We adopt the notion that
$$
|\xi_A|\geq |\xi_B|\geq |\xi_C|\geq |\xi_D|\geq |\xi_E|\geq |\xi_F|,
$$
and let
\begin{equation*}\aligned
\Omega_1:=&\;\big\{(\xi_1,\cdots, \xi_6)\in \Gamma_6:  |\xi_A|\sim
|\xi_B|\gtrsim N\gg |\xi_C|,|\xi_A^3+\xi_B^3|\gg
|\xi_C^3+\cdots+\xi_F^3|\big\};\\
\Omega_2:=&\;\big\{(\xi_1,\cdots, \xi_6)\in \Gamma_6:
|\xi_C|\gtrsim N,|\xi_C|\gg |\xi_D|\big\};\\
\Omega_3:=&\;\big\{(\xi_1,\cdots, \xi_6)\in \Gamma_6:
 |\xi_D|\gtrsim N\gg |\xi_E|,|\xi_A+\xi_B|\gg |\xi_E+\xi_F|,\\
&\qquad\qquad\qquad\qquad\;\;
|m^2(\xi_A)\xi_A^3+\cdots+m^2(\xi_D)\xi_D^3|\gg|\xi_E^3+\xi_F^3|,\\
&\qquad\qquad\qquad\qquad\;\;
|\xi_A+\xi_B||\xi_A+\xi_C||\xi_B+\xi_C|\gg |\xi_E||\xi_A|^2\big\}.
\endaligned\end{equation*}
Then we rewrite (\ref{dE1}) by
$$
\dfrac{d}{dt}E^1_I(u(t))=\Lambda_{6}(\bar{M}_6)
+\Lambda_{6}(\tilde{M}_6)+\Lambda_{10}(M_{10}),
$$
where $\bar{M}_6$ and $\tilde{M}_6$ are defined by
\begin{equation}\label{eq-1}
\begin{split}
\bar{M}_6 &:= (\chi_{\Gamma_6} -\chi_\Omega )M_6^1;\\
\tilde{M}_6 &:= \chi_\Omega M_6^1+\chi_{\Gamma_6} M_6^2;\\
\Omega &:=\Omega_1\cup\Omega_2\cup\Omega_3.
\end{split}
\end{equation}
Now we are ready to define a new modified energy $E^2_I(u(t))$ by
\begin{equation}
E^2_I(u(t)):=\Lambda_6(\tilde{\sigma}_6)+E^1_I(u(t)), \label{E2}
\end{equation}
where
\begin{equation}\label{eq-2}
\tilde{\sigma}_6:=-\tilde{M}_6/\alpha_6.
\end{equation}
Then by applying the differentiation formula \eqref{dLMD} again, we see that
\begin{equation}
\dfrac{d}{dt}E^2_I(u(t))=\Lambda_{6}(\bar{M}_6)+\Lambda_{10}(\bar{M}_{10}),\label{dE2}
\end{equation}
where
\begin{equation}\label{eq-3}
\begin{split}
\bar{M}_{10}:=&-6i\big[\big(\sigma_6(\xi_1,\cdots,\xi_{5},\xi_{6}+\cdots+\xi_{10})+\\
&\quad +\tilde{\sigma}_6(\xi_1,\cdots,\xi_{5},\xi_{6}+\cdots+\xi_{10}))(\xi_{6}+\cdots+\xi_{10}\big)\big]_{sym}.
\end{split}
\end{equation}
\begin{remark}
On the support of ``resonant" $\bar{M}_6$, when $|\xi_B|\gtrsim N\gg |\xi_C|$, one may find that $\bar{M}_6$ is controlled by $|\xi_C^3+\cdots+\xi_F^3|\lesssim |\xi_C||\xi_D||\xi_E|$ (See estimates of $A_1$ in the proof of Proposition \ref{6-linear}), which is a relative low frequency term and would be expected to give a better decay due to the high/low interaction. On the other hand, on the support of $\tilde{M}_6$, as is shown in Lemma \ref{le:6-linear-bound}, $$|\tilde{M}_6|\lesssim|\alpha_6|,$$ which is non-resonant and will give a small error in energy-increment. We will show that $E^2_I(u(t))$ is almost conserved from the multilinear estimates, which in turn shows that $E^1_I(u(t))$ is almost conserved by ignoring an error coming from $\tilde{\sigma}_6$.
\end{remark}

To prove Theorem \ref{thm:main}, it suffices to prove the following theorem.
\begin{thm}[Existence of an almost conserved quantity]\label{thm:main-2}
Let the notation be as above. Then for a solution $u$ to \eqref{gkdv} which is smooth-in-time, Schwarz-in-space on a time interval $[0,\delta]$, we have
\begin{itemize}
\item (Fixed-time bound) For $1/3<s<1$,
\begin{equation}
\left| \Lambda_6(\tilde{\sigma}_6)(t)\right| \lesssim N^{0-}\|Iu(t)\|^6_{H^1_x}.
\end{equation}
\item (Almost conservation law) For $t\in [0,\delta]$,
 \begin{equation}
 \left|E^2_I(u(t))-E^2_I(u_0)\right|\le CN^{-\frac 72+}\bigl(\|Iu\|^{6}_{X^\delta_{1,\frac 12+}}+\|Iu\|^{10}_{X^\delta_{1,\frac 12+}}\bigr).
 \end{equation}
\end{itemize}
\end{thm}

We will show that Theorem \ref{thm:main-2} implies Theorem \ref{thm:main} in Section \ref{sec:final-argument}. Now we focus on establishing the claims in Theorem \ref{thm:main-2}, which will occupy the next section.

\section{Fixed-time bound and almost conservation law}\label{sec:main-arguments}
In this section, we prove Theorem \ref{thm:main-2}. We start with a few basic facts which will be only used in this section and are taken from \cite{CKSTT-03-KDV}. The first is the following well-known arithmetic fact \cite[(4.2)]{CKSTT-03-KDV}.
\begin{equation}\label{eq:arithmetic}
\xi_1+\xi_2+\xi_3+\xi_4=0 \Rightarrow \xi_1^3+\xi_2^3+\xi_3^3+\xi_4^3=3(\xi_1+\xi_2)(\xi_1+\xi_3)(\xi_1+\xi_4).
\end{equation}
Then we record the following forms of the mean value theorem. To prepare for it, we state a definition: Let $a$ and $b$ be two smooth functions of real variables. We say that $a$ is controlled by $b$ if $b$ is non-negative and satisfies $b(\xi)\sim b(\xi')$ for $|\xi|\sim |\xi'|$ and $$a(\xi) \lesssim b(\xi), \,a^\prime (\xi) \lesssim \frac {b(\xi)}{|\xi|},\,a^{\prime\prime} \lesssim \frac {b(\xi)}{|\xi|^2}.$$
\begin{lem}\label{le:mvt}
If $a$ is controlled by $b$ and $|\eta|, |\lambda|\ll |\xi|$, then
\begin{itemize}
\item (Mean value theorem)
\begin{equation}
\left| a(\xi+\eta)-a(\xi)\right| \lesssim |\eta| \frac {b(\xi)}{|\xi|}.
\end{equation}
\item (Double mean value theorem)
\begin{equation}
\left| a(\xi+\eta+\lambda)-a(\xi+\eta)-a(\xi+\lambda)+a(\xi)\right| \lesssim |\eta||\lambda| \frac {b(\xi)}{|\xi|^2}.
\end{equation}
\end{itemize}
\end{lem}
We will use this lemma in a context that $a(\xi)=m^2(\xi)|\xi|^3$
for $m$  defined above with the choice $b(\xi)=4m^2(\xi)|\xi|^3$.

\subsection{Fixed-time bound}
The first part of Theorem \ref{thm:main-2} is a consequence of the following two lemmas.
\begin{lem}\label{le:6-linear-bound}
Let the notations be as in \eqref{eq-1} and \eqref{eq-2}. Then
\begin{equation}\label{eq-4}
|\tilde{M}_6|\lesssim|\alpha_6|, \text{ i.e., }
|\tilde{\sigma}_6|\lesssim 1.
\end{equation}
\end{lem}
\begin{proof}
By a symmetry consideration, we may assume that
\begin{equation*}
\aligned \big|\xi_{1}\big|\geq   \big|\xi_{2} \big| \geq \big|
\xi_{3} \big| \geq  \big|\xi_{4} \big| \geq  \big|\xi_{5} \big| \geq
\big|\xi_{6} \big|.
\endaligned
\end{equation*}
By the definitions of $\tilde{M}_6$ and $M_6^2$, we only need to
show
\begin{equation}\label{eq-6}
\left|\chi_\Omega M_6^1\right|\lesssim |\alpha_6|,
\end{equation}since $|\sigma_6|\lesssim 1$ always holds.
We will prove this bound case by case by analyzing it on domains $\Omega_i$ for $i=1,2,3$.  Recall that
\begin{equation*}
\begin{split}
|M_6^1| & \sim |m^2(\xi_1)\xi_1^3+\cdots +m^2(\xi_6)\xi_6^3|.\\
|\alpha_6| &\sim |\xi_1^3+\cdots +\xi_6^3|.\\
\text{On } \Gamma_6, \quad&\xi_1+\cdots+\xi_6=0.
\end{split}
\end{equation*}

\textbf{Case 1.} On $\Omega_1$, there holds that $|\xi_3^3+\cdots+\xi_6^3|\ll |\xi_1^3+\xi_2^3|$ and $m(\xi_3)=\cdots=m(\xi_6)\equiv 1$. It follows that
\begin{equation}\label{eq-7}
|\alpha_6|\gtrsim |\xi_1^3+\xi_2^3|\gtrsim |\xi_1+\xi_2|(\xi_1^2+\xi_2^2),
\end{equation} as $|\xi_1\xi_2|\le (\xi_2^2+\xi_2^2)/2$ always holds and then $\xi_1^2+\xi^2-\xi_1\xi_2\ge (\xi_2^2+\xi_2^2)/2$.

On the other hand, by the mean value theorem in Lemma \ref{le:mvt}, and the fact that $m$ is even and $|\xi_1-(-\xi_2)|\lesssim |\xi_3|\ll |\xi_i|$ for $i=1,2$,
\begin{equation}\label{eq-8}
\aligned
|M_6^1| \sim &\;
\left|m^2(\xi_1)\xi_1^3+m^2(\xi_2)\xi_2^3+\xi_3^3+\cdots+\xi_6^3\right|\\
\lesssim &\;\left|m^2(\xi_1)\xi_1^3-m^2(-\xi_2)(-\xi_2)^3\right|+\left|\xi_3^3+\cdots+\xi_6^3\right|\\
\lesssim &\;
\left|(\xi_1+\xi_2)m^2(\xi_2)|\xi_2|^2\right|+\left|\xi_3^3+\cdots+\xi_6^3\right|\\
\lesssim&\;
|\xi_1+\xi_2|(\xi_1^2+\xi_2^2)+|\xi_3^3+\cdots+\xi_6^3|\\
\lesssim&\; |\xi_1+\xi_2|(\xi_1^2+\xi_2^2).
\endaligned\end{equation}Thus \eqref{eq-6} follows from \eqref{eq-7} and \eqref{eq-8}.

\textbf{Case 2.} On $\Omega_2$, since $|\xi_3|\gg |\xi_4|$, there always holds that $|\xi_1+\xi_2|\sim |\xi_3|$; moreover we have $\xi_1\cdot\xi_2<0$; otherwise, if $\xi_1$ and $\xi_2$ would have the same sign, then from the support information of $\Gamma_6$ and $|\xi_4|\ll |\xi_3|$,
$$2|\xi_2|\le |\xi_1+\xi_2|=|\xi_3+\xi_4+\xi_5+\xi_6|\le 3|\xi_3|/2\le 3|\xi_2|/2,$$
 which is obviously a contradiction since $|\xi_2|>0$. Then
\begin{equation}\label{eq-9}
\aligned
|\alpha_6| \gtrsim&\;
\left|(\xi_1+\xi_2)(\xi_1^2+\xi_2^2-\xi_1\xi_2)+(\xi_3+\xi_4)(\xi_3^2+\xi_4^2-\xi_3\xi_4)\right|
+o(|\xi_3|(\xi_1^2+\xi_2^2))\\
=&\;
\left|(\xi_1+\xi_2)(\xi_1^2+\xi_2^2-\xi_1\xi_2)+\bigl(-(\xi_1+\xi_2)-(\xi_5+\xi_6)\bigr)(\xi_3^2+\xi_4^2-\xi_3\xi_4)\right|
+o(|\xi_3|(\xi_1^2+\xi_2^2))\\
\gtrsim&\;
\left|(\xi_1+\xi_2)(\xi_1^2+\xi_2^2-\xi_1\xi_2-3\xi_3^2/2)\right|+o(|\xi_3|(\xi_1^2+\xi_2^2))\\
= &\;
\left|(\xi_1+\xi_2)\bigl((\xi_1^2+\xi_2^2-\xi_3^2/2)+(-\xi_1\xi_2-\xi_3^2)\bigr)\right|+o(|\xi_3|(\xi_1^2+\xi_2^2))\\
\gtrsim &\; |\xi_3|(\xi_1^2+\xi_2^2).
\endaligned\end{equation} The first inequality follows since $|\xi_5^3+\xi_6^3|\ll |\xi_3|(\xi_1^2+\xi_2^2)$,
 the third follows since $|\xi_3^2+\xi_4^2-\xi_3\xi_4|\le 3|\xi_3|^2/2$ and $\left|(\xi_5+\xi_6)(\xi_3^2+\xi_4^2-\xi_3\xi_4)\right|\ll |\xi_3|(\xi_1^2+\xi_2^2)$,
  and the last follows since $(\xi_1^2+\xi_2^2-\xi_3^2/2)\ge (\xi_1^2+\xi_2^2)/2$ and $-\xi_1\xi_2-\xi_3^2>0$.

On the other hand, by the same use of mean value theorem as in \eqref{eq-8}, we see that
\begin{equation}\label{eq-10}
\aligned
|M_6^1|\leq &\;
\left|m^2(\xi_1)\xi_1^3+m^2(\xi_2)\xi_2^3\right|+|m^2(\xi_3)\xi_3^3|+|\xi_4^3+\xi_5^3+\xi_6^3|\\
\lesssim&\;|\xi_1+\xi_2||\xi_2^2|+|\xi_3|^3\\
\lesssim&\; |\xi_3|(\xi_1^2+\xi_2^2),
\endaligned\end{equation} since $|\xi_1+\xi_2|\lesssim |\xi_3|$. Thus \eqref{eq-6} follows from \eqref{eq-9} and \eqref{eq-10} again.

\textbf{Case 3.} We consider $\Omega_3$. In fact we will consider $\Omega_3\setminus\Omega_2$ by \textbf{Case 2}. We split $\Omega_3\setminus\Omega_2$ into two parts:
\begin{equation}\label{O3}
\textbf{Subcase 3a}: |\xi_1|\gg|\xi_3|\sim|\xi_4|;\qquad \textbf{Subcase 3b}:
|\xi_1|\sim|\xi_2|\sim|\xi_3|\sim|\xi_4|.
\end{equation}

\textbf{Subcase 3a}. On $\Omega_3$, there holds that $|\xi_1+\xi_2|\gg |\xi_5+\xi_6|$. Thus by the same consideration as in \eqref{eq-9},
\begin{equation}\label{eq-11}
\begin{split}
|\alpha_6|&\gtrsim \left|\xi_1^3+\xi_2^3+\xi_3^3+\xi_4^3\right|+o(|\xi_1+\xi_2|\xi_1^2)\\
&=\left|(\xi_1+\xi_2)(\xi_1^2+\xi_2^2-\xi_1\xi_2)+\bigl(-(\xi_1+\xi_2)-(\xi_5+\xi_6)\bigr)
(\xi_3^2+\xi_4^2-\xi_3\xi^4)\right|+o(|\xi_1+\xi_2|\xi_1^2)\\
&\gtrsim \left|(\xi_1+\xi_2)(\xi_1^2+\xi_2^2-\xi_1\xi_2-\xi^2_3-\xi_4^2+\xi_3\xi_4)\right|+o(|\xi_1+\xi_2|\xi_1^2)\\
&\gtrsim |\xi_1+\xi_2|\xi_1^2.
\end{split}
\end{equation} The last inequality follows since $\xi_1^2+\xi_2^2-\xi_1\xi_2 \ge (\xi_1^2+\xi_2^2)/2$ and $|\xi_4|\le |\xi_3|\ll |\xi_1|\sim |\xi_2|.$

On the other hand, we have
\begin{equation}\label{eq-12}
\begin{split}
|M_6^1|&\lesssim \left|m^2(\xi_1)\xi_1^3+m^2(\xi_2)\xi_2^3\right| +\left|m^2(\xi_3)\xi_3^3+m^2(\xi_4)\xi_4^3\right|+ |\xi_5^3+\xi_6^3|\\
&\lesssim |\xi_1+\xi_2|\xi_2^2+|\xi_3+\xi_4|\xi_3^2+|\xi_5+\xi_6|\xi_5^2\\
&\lesssim |\xi_1+\xi_2|\xi_1^2,
\end{split}
\end{equation} by the mean value theorem in Lemma \ref{le:mvt} and the fact that $|\xi_5+\xi_6|\ll |\xi_1+\xi_2|\sim |\xi_3+\xi_4|$ on $\Omega_3$.

\textbf{Case 3b}. We may assume that $\xi_1> 0$ by symmetry. Then one of the following four subcases always occurs:
\begin{enumerate}
\item[3b-I:] $ \xi_1>0, \xi_2<0, \xi_3<0, \xi_4<0;$
\item[3b-II:] $ \xi_1>0, \xi_2<0, \xi_3>0, \xi_4<0;$
\item[3b-III:] $ \xi_1>0, \xi_2<0, \xi_3<0, \xi_4>0;$
\item[3b-IV:] $ \xi_1>0, \xi_2>0, \xi_3<0, \xi_4<0.$
\end{enumerate}
In fact, as $|\xi_2|\le |\xi_1|$, $\xi_1>0$ will imply that $\xi_1+\xi_2\ge 0$. On $\Omega_3$, there holds that
$$|\xi_1+\xi_2+\xi_3+\xi_4|=|\xi_5+\xi_6|\ll |\xi_1+\xi_2|\sim |\xi_3+\xi_4|,$$
which implies that $(\xi_1+\xi_2)$ and $(\xi_3+\xi_4)$ has different signs, i.e., $\xi_3+\xi_4\le 0$. This information, together with that $|\xi_3|\le |\xi_2|\le |\xi_1|$, implies the classification above.

\textbf{Case 3b-I}. $\xi_1\xi_2<0$ and $\xi_3\xi_4>0$. Since $|\xi_5+\xi_6|\ll |\xi_1+\xi_2|\sim |\xi_3+\xi_4|$,
by the same consideration as in \eqref{eq-11}, we have
\begin{equation}\label{eq-14}
\begin{split}
|\alpha_6|&\gtrsim |\xi_1+\xi_2||\xi_1^2+\xi_2^2-\xi_1\xi_2-\xi_3^2-\xi_4^2+\xi_3\xi_4|+o(|\xi_1+\xi_2|\xi_1^2)\\
&\gtrsim |\xi_1+\xi_2|\xi_1^2
\end{split}
\end{equation} as $-\xi_1\xi_2-\xi_3^2\ge 0$ and $\xi_3\xi_4-\xi_4^2\ge 0$. On the other hand, the usual mean value theorem as in \eqref{eq-12} will give that
\begin{equation}\label{eq-15}
|M_6^1|\lesssim |\xi_1+\xi_2|\xi_1^2.
\end{equation}
Thus \eqref{eq-6} follows again from \eqref{eq-14} and \eqref{eq-15}.

For \textbf{Case 3b-II} to \textbf{Case 3b-IV}, we denote $\overline{\xi_4}=\xi_4+\xi_5+\xi_6$, and
rewrite
\begin{equation}\label{eq-16}
\begin{split}
m^2(\xi_1)\xi_1^3+\cdots+m^2(\xi_4)\xi_4^3 &=
m^2(\xi_1)\xi_1^3+m^2(\xi_2)\xi_2^3+m^2(\xi_3)\xi_3^3\\
&\quad +m^2(\overline{\xi_4})\overline{\xi_4}^3-\left(m^2(\overline{\xi_4})
\overline{\xi_4}^3-m^2(\xi_4)\xi_4^3\right).
\end{split}
\end{equation}
On one hand, by the double mean value theorem in Lemma \ref{le:mvt}, we have
\begin{equation}\label{eq-17}
\left|m^2(\xi_1)\xi_1^3+m^2(\xi_2)\xi_2^3+m^2(\xi_3)\xi_3^3
+m^2(\overline{\xi_4})\overline{\xi_4}^3\right|\lesssim
m^2(\xi_1)|(\xi_1+\xi_2)(\xi_1+\xi_3)(\xi_2+\xi_3)|
\end{equation}
Indeed, suppose we are in \textbf{Case 3b-II}, $\xi_1, \xi_3>0$ and $\xi_2, \xi_4<0$; then by taking $\xi=\xi_1$,
$\eta=-\xi_1-\xi_2$, $\lambda=\xi_2+\xi_3$, we apply the double mean value theorem in Lemma \ref{le:mvt},
\begin{equation*}
\left|m^2(\xi_1)\xi_1^3+m^2(\xi_2)\xi_2^3+m^2(\xi_3)\xi_3^3
+m^2(\overline{\xi_4})\overline{\xi_4}^3\right|\lesssim
m(\xi_1)|\xi_1||(\xi_1+\xi_2)(\xi_2+\xi_3)|.
\end{equation*} Then by using the fact that $\xi_1,\xi_3$ have the same signs, we see that  $|\xi_1|\le |\xi_1+\xi_3|$.  So there follows \eqref{eq-17}. The other two cases are treated similarly.

On the other hand, by the mean value theorem, we obtain
\begin{equation}\label{eq-18}
\left|m^2(\overline{\xi_4})\overline{\xi_4}^3-m^2(\xi_4)\xi_4^3\right|
\lesssim m^2(\xi_4)|\xi_5+\xi_6||\xi_4|^2\lesssim
m^2(\xi_4)|\xi_5||\xi_4|^2.
\end{equation}
Since $|\xi_1+\xi_2||\xi_1+\xi_3||\xi_2+\xi_3|\gg |\xi_5||\xi_1|^2$
and $|\xi_1|\sim |\xi_4|$, by (\ref{eq-16}), we have
\begin{equation}\label{eq-19}
|m^2(\xi_1)\xi_1^3+\cdots+m^2(\xi_4)\xi_4^3|\lesssim
m^2(\xi_1)|(\xi_1+\xi_2)(\xi_1+\xi_3)(\xi_2+\xi_3)|.
\end{equation}
By the definition of $M_6^1$ and the support information of $\Omega_3$, it follows that
\begin{align}
& \label{eq-20} |M_6^1| \lesssim m(\xi_1)|(\xi_1+\xi_2)(\xi_1+\xi_3)(\xi_2+\xi_3)|\le |(\xi_1+\xi_2)(\xi_1+\xi_3)(\xi_2+\xi_3)|,\\
&\label{eq-21}
|\xi_5^3+\xi_6^3| \ll m^2(\xi_1)|(\xi_1+\xi_2)(\xi_1+\xi_3)(\xi_2+\xi_3)|\le |(\xi_1+\xi_2)(\xi_1+\xi_3)(\xi_2+\xi_3)|.
\end{align}

Reasoning similarly as in proving \eqref{eq-19}, we have
\begin{equation}\label{eq-22}
|\alpha_6|\gtrsim |(\xi_1+\xi_2)(\xi_1+\xi_3)(\xi_2+\xi_3)|.
\end{equation}
Indeed, let $\bar{\xi}_4$ be defined as above, then
$\alpha_6\ge |\xi_1^3+\xi_2^3+\xi_3^3+\bar{\xi_4}^3|-|\bar{\xi}^3_4-\xi_4^3|-|\xi^3_5+\xi^3_6|$.
Then by the definition of $\Omega_3$ and \eqref{eq-21}, $$|\bar{\xi}^3_4-\xi_4^3|, |\xi^3_5+\xi^3_6|\ll |(\xi_1+\xi_2)(\xi_1+\xi_3)(\xi_2+\xi_3)|,$$
while $|\xi_1^3+\xi_2^3+\xi_3^3+\bar{\xi_4}^3|\sim |(\xi_1+\xi_2)(\xi_1+\xi_3)(\xi_2+\xi_3)|$ by the arithmetical fact \eqref{eq:arithmetic}.

Together with \eqref{eq-20}, \eqref{eq-22} will yield
\begin{equation}\label{eq-51}
|\alpha_6|\gtrsim |(\xi_1+\xi_2)(\xi_1+\xi_3)(\xi_2+\xi_3)|\gtrsim |M_6^1|.
\end{equation}
Hence \eqref{eq-6} follows again. This completes the proof for Lemma \eqref{le:6-linear-bound}.
\end{proof}
Now we can establish the first part of Theorem \ref{thm:main-2}.
\begin{lem}\label{le:pointwise-bound}
For any $1>s>1/3$, we have
\begin{equation}\label{eq-23}
|\Lambda_6(\tilde{\sigma}_6)(t)|\lesssim N^{0-}\>\|Iu(t)\|_{H^1_x}^6.
\end{equation}
\end{lem}
\begin{proof}
By Lemma \ref{le:6-linear-bound}, it suffices to show
\begin{equation}\label{3.9}
\displaystyle\int_{\Gamma_6} \frac{\left|\widehat{f_1}(\xi_1,t)\cdots \widehat{
f_6}(\xi_6,t)\right|} {\langle\xi_1\rangle m(\xi_1)\cdots
\langle\xi_6\rangle m(\xi_6)} \lesssim
N^{0-}\>\|f_1(t)\|_{L^2_{x}}\cdots \|f_6(t)\|_{L^2_{x}}.
\end{equation}
By the definition of $m$, we have
$$
\langle\xi\rangle m(\xi)= \langle\xi\rangle,\mbox{ for }|\xi|\leq
N;\quad \langle\xi\rangle m(\xi)\sim N^{1-s}|\xi|^s, \mbox{ for
}|\xi|\gtrsim N,
$$
By using Remark \ref{re:two-reductions}, we may assume that $|\xi_1|\ge \cdots \ge |\xi_6|$ and $|\xi_1|\sim |\xi_2|\gtrsim N$. To estimate the left hand side of (\ref{3.9}), we may assume that the spatial Fourier transforms of the $f_i$ are nonnegative and consider the worst case where $|\xi_6|\gtrsim N$, by using Plancherel's theorem in the spatial variable,
\begin{equation}\label{eq-24}
\begin{split}
&\lesssim N^{6(s-1)} \int \frac{\widehat{f_1}(\xi_1)\cdots
\widehat{f_6}(\xi_6)} {\langle\xi_1\rangle^s \cdots
\langle\xi_6\rangle^s} \\
&\lesssim N^{6(s-1)} \int \frac{\widehat{f_1}(\xi_1)\cdots \widehat{f_6}(\xi_6)} {\langle\xi_3\rangle^{3s/2} \cdots\langle\xi_6\rangle^{3s/2}}\\
&\lesssim N^{0-}\|f_1(t)\|_{L^2_{x}}\|f_2(t)\|_{L^2_{x}}
\Big\|J^{-\frac{1}{2}-}_xf_3(t)\Big\|_{L^\infty_{x}}\cdots
\Big\|J^{-\frac{1}{2}-}_xf_6(t)\Big\|_{L^\infty_{x}},
\end{split}
\end{equation} for $1/3<s<1$. Hence \eqref{eq-23} follows from Sobolev's inequality.
\end{proof}

\subsection{An upper bound on the increment of $E_I^2(u(t))$} In the subsection, we will establish the second half of Theorem \ref{thm:main-2}. By the multilinear correction analysis, the almost conservation law
of $E_I^2(u(t))$ is the key ingredient in the proof of the global well-posedness below the energy space. While by (\ref{dE2}), the main process to estimate the increment of $E_I^2(u(t))$ is the following 6-linear and 10-linear estimates.

\begin{prop}\label{6-linear}
For any $s>3/8$, we have
\begin{equation}
 \left|\displaystyle\int_0^\delta\Lambda_6(\bar{M}_6)\,dt\right|
 \lesssim N^{-\frac{7}{2}+}\>\|Iu\|_{X_{1,\frac{1}{2}+}^\delta}^6.
\label{eq-52}
\end{equation}
\end{prop}
\begin{proof} By using Remark \ref{re:two-reductions}, we may assume that
\begin{equation*}
\begin{split}
&\big|\xi_{1}\big|\geq   \big|\xi_{2} \big| \geq \big|
\xi_{3} \big| \geq  \big|\xi_{4} \big| \geq  \big|\xi_{5} \big| \geq
\big|\xi_{6} \big|,\\
&|\xi_i|\sim N_i, i=1,\cdots, 6,\text{ and } N_1\sim N_2\gtrsim N.
\end{split}
\end{equation*}
To recover the sum at the end we need to borrow a $N_1^{0-}$ but this will not be mentioned and it will only be recorded at the end by paying a price equivalent $N^{0+}$. By Plancherel's theorem, we only need to show
\begin{equation*}
\left|\displaystyle\int_0^\delta\int_{\Gamma_6}
\frac{\bar{M}_6(\xi_1,\cdots,\xi_6)\widehat{f_1}(\xi_1)\cdots \widehat{f_6}(\xi_6)}
{\langle\xi_1\rangle m(\xi_1)\cdots \langle\xi_6\rangle m(\xi_6)}\right|
\lesssim N^{-\frac{7}{2}+}\>\|f_1\|_{X^\delta_{0,\frac{1}{2}+}}\cdots
\|f_6\|_{X^\delta_{0,\frac{1}{2}+}}.
\end{equation*} Because of the definition $\|u\|_{X^\delta_{s,b}([0,\delta]\times \mathbb{R})}:=\inf\{\|U\|_{X_{s,b}(\R\times \R)}:\,U|_{[0,\delta]}=u\}$, it reduces to show that
\begin{equation}\label{3.11}
\displaystyle\int dt\int_{\Gamma_6}
\frac{\left|\bar{M}_6(\xi_1,\cdots,\xi_6)\widehat{f_1}(\xi_1)\cdots \widehat{f_6}(\xi_6)\right|}
{\langle\xi_1\rangle m(\xi_1)\cdots \langle\xi_6\rangle m(\xi_6)}
\lesssim N^{-\frac{7}{2}+}\>\|f_1\|_{X_{0,\frac{1}{2}+}}\cdots
\|f_6\|_{X_{0,\frac{1}{2}+}}.
\end{equation} In view of this inequality, we may assume that the spatial Fourier transforms of the $f_i$ are nonnegative, which will be used in the arguments throughout the paper without being mentioned. Now we split it into five regions:
\begin{equation*}\aligned
  A_1 :=& \{\xi\in (\Gamma_6\setminus\Omega): |\xi_2|\gtrsim N\gg |\xi_3|\}; \\
  A_2 :=& \{\xi\in (\Gamma_6\setminus\Omega): |\xi_3|\gtrsim N\gg |\xi_4|\}; \\
  A_3 :=& \{\xi\in (\Gamma_6\setminus\Omega): |\xi_4|\gtrsim N\gg |\xi_5|\}; \\
  A_4 :=& \{\xi\in (\Gamma_6\setminus\Omega): |\xi_5|\gtrsim N\gg |\xi_6|\}; \\
  A_5 :=& \{\xi\in (\Gamma_6\setminus\Omega): |\xi_6|\gtrsim N\}.
\endaligned\end{equation*}

\textbf{Estimate in $A_1$.} In $A_1$, there holds that
\begin{equation}\label{eq-25}
|\xi_1^3+\xi_2^3|\lesssim |\xi_3^3+\cdots+\xi_6^3|;
\end{equation} also $|\xi_1+\xi_2|=|\xi_3+\cdots +\xi_6|$ and $|\xi_3|\ll |\xi_2|$ imply that $\xi_1\xi_2<0$.
The two estimates above imply that
\begin{equation}\label{eq-26}
|\xi_1+\xi_2|=\frac {|\xi_1^3+\xi_2^3|}{|\xi_1^2+\xi_2^2-\xi_1\xi_2|} \lesssim |\xi_3|\frac {\xi_3^2}{\xi_1^2+\xi_2^2}\ll |\xi_3|.
\end{equation}
Then by the mean value theorem in Lemma \ref{le:mvt}, we have
\begin{equation*}\aligned
|\bar{M}_6|&\leq \big|m^2(\xi_1)\xi_1^3+m^2(\xi_2)\xi_2^3\big|+\big|\xi_3^3+\cdots +\xi_6^3\big|\\
&\lesssim \big|\xi_1+\xi_2\big|\xi_2^2+\big|\xi_3^3+\cdots +\xi_6^3\big|\\
&\lesssim \big|\xi_1+\xi_2\big|\bigl|\xi_1^2+\xi_2^2-\xi_1\xi_2\bigr|+\big|\xi_3^3+\cdots +\xi_6^3\big|\\
&\lesssim \big|\xi_3^3+\cdots+\xi_6^3\big|\lesssim
|\xi_3||\xi_4||\xi_5|.
\endaligned\end{equation*}
Therefore, by Plancherel's theorem in the spatial variable, H\"older's inequality followed by Lemma \ref{le:strichartz-esimate-2} and Lemma \ref{le:strichartz-esimate-3}, the left-hand side of (\ref{3.11}) is bounded by
\begin{equation}\label{eq-27}
\aligned
&\lesssim N^{2s-2}\int dt\int_{A_1}
\frac{\widehat{f_1}(\xi_1)\cdots \widehat{f_6}(\xi_6)}
{|\xi_1|^s |\xi_2|^s\langle\xi_6\rangle}\\
\lesssim &\; N^{2s-2} \int dt \int_{\Gamma_6}|\xi_1|^{-s-1}|\xi_2|^{-s-1}\langle \xi_5\rangle^{\frac 12+} \left(|\xi_1|\widehat{f_1}\widehat{f_3}\right)
\left( |\xi_2|\widehat{f_2}\widehat{f_4}\right)
\left(\langle \xi_5\rangle^{-\frac12-}\widehat{f_5}\right)\left(\langle \xi_6\rangle^{-1}\widehat{f_6}\right)\\
\lesssim &\; N^{-\frac 72+} \int dt \int_{\Gamma_6}\left(|\xi_1|\widehat{f_1}\widehat{f_3}\right)
\left( |\xi_2|\widehat{f_2}\widehat{f_4}\right)
\left(\langle \xi_5\rangle^{-\frac12-}\widehat{f_5}\right)\left(\langle \xi_6\rangle^{-1}\widehat{f_6}\right)\\
\lesssim& \;
N^{-\frac{7}{2}+}\int dt\int \bigl[(D_xf_1)f_3\bigr]\bigl[(D_xf_2)f_4\bigr]\bigl[J_x^{-\frac 12-}f_5\bigr] \bigl[J_x^{-\frac 12-}f_6\bigr] dx\\
\lesssim& \;
N^{-\frac{7}{2}+}\left\|D_xf_1f_3\right\|_{L^2_{x,t}}
\left\|D_xf_2f_4\right\|_{L^2_{x,t}}
\left\|J_x^{-\frac{1}{2}-}f_5\right\|_{L^\infty_{x,t}}
\left\|J_x^{-1}f_6\right\|_{L^\infty_{x,t}}\\
\lesssim&\; N^{-\frac{7}{2}+}\> \|f_1\|_{X_{0,\frac{1}{2}+}}\cdots
\|f_6\|_{X_{0,\frac{1}{2}+}}, \endaligned\end{equation}
 where we have used the fact that $|\xi_3|\ll |\xi_1|$ and $|\xi_4|\ll |\xi_2|$.

\textbf{Estimate in $A_2$}. Note that $A_2=\emptyset$,
thus $\bar{M}_6=0$.

\textbf{Estimate in $A_3$}. We may split $A_3$ into
three regions again,
\begin{equation*}\aligned
  A_{31} :=&\; \{\xi\in A_3:\,|\xi_1+\xi_2|\lesssim |\xi_5+\xi_6|\}; \\
  A_{32} :=&\; \{\xi\in A_3:\,|m^2(\xi_1)\xi_1^3+\cdots+m^2(\xi_4)\xi_4^3|\lesssim|\xi_5^3+\xi_6^3|\}; \\
  A_{33} :=&\; \{\xi\in A_3\setminus(A_{31}\cup A_{32}):\,|(\xi_1+\xi_2)(\xi_1+\xi_3)(\xi_2+\xi_3)|\lesssim
|\xi_5||\xi_1|^2\}. \endaligned
\end{equation*} We claim that in all three cases,
\begin{equation}\label{3.17}
|\bar{M}_6|\lesssim |\xi_5||\xi_1^2+\xi_2^2|.
\end{equation}
Indeed, \begin{itemize}
\item In $A_{31}$, we also have $|\xi_3+\xi_4|\lesssim |\xi_5+\xi_6|$, which implies that
\begin{equation*}\aligned
|\bar{M}_6| \leq &\;
|m^2(\xi_1)\xi_1^3+m^2(\xi_2)\xi_2^3|+|m^2(\xi_3)\xi_3^3+m^2(\xi_4)\xi_4^3|
+|\xi_5^3+\xi_6^3|\\
\lesssim&\;
|\xi_1+\xi_2|(\xi_1^2+\xi_2^2)+|\xi_3+\xi_4|(\xi_3^2+\xi_4^2)
+|\xi_5+\xi_6|(\xi_5^2+\xi_6^2)\\
\lesssim&\; |\xi_5|(\xi_1^2+\xi_2^2).
\endaligned\end{equation*}
\item In $A_{32}$, $|\bar{M}_6|\lesssim |\xi_5^3+\xi_6^3|\lesssim |\xi_5|(\xi_1^2+\xi_2^2).$
\item In $A_{33}$, since $|\xi_3|\sim |\xi_4|$ in $\Gamma\setminus\Omega_2$,
we split it into two case as in (\ref{O3}).
\begin{itemize}
\item If $|\xi_1|\gg |\xi_3|\sim |\xi_4|$, then $|\xi_1+\xi_3|, |\xi_1+\xi_4|\sim |\xi_1|$. Also the same reasoning as proving (\ref{eq-12}) implies that
$$
|\bar{M}_6|\lesssim|\xi_1+\xi_2||\xi_1|^2 \sim
|(\xi_1+\xi_2)(\xi_1+\xi_3)(\xi_2+\xi_3)|\lesssim |\xi_5||\xi_1|^2.
$$
The last inequality follows from the support of $A_{33}$.
\item If $|\xi_1|\sim |\xi_4|$, note that \eqref{eq-16} to \eqref{eq-19} still hold in this
situation; hence (\ref{3.17}) follows because $|(\xi_1+\xi_2)(\xi_1+\xi_3)(\xi_2+\xi_3)|\lesssim|\xi_5||\xi_1|^2$.
\end{itemize}
\end{itemize}
Therefore, by (\ref{3.17}) and Lemmas \ref{le:strichartz-esimate-1}, \ref{le:strichartz-esimate-2} and \ref{le:strichartz-esimate-3}, for $s>3/8$, by a similar reasoning as in \eqref{eq-27},
\begin{equation*}\aligned
\text{LHS of }\eqref{3.11} \lesssim&\; N^{4s-4}\displaystyle\int dt \int_{A_3}
\frac{\widehat{f_1}(\xi_1)\cdots \widehat{f_6}(\xi_6)}
{|\xi_1|^{s-1} |\xi_2|^{s-1}|\xi_3|^{s} |\xi_4|^{s}\langle\xi_6\rangle}\\
\lesssim&\;
N^{-\frac{7}{2}}\left\|D_xf_1f_5\right\|_{L^2_{x,t}}
\left\|D_xf_2\right\|_{L^\infty_{x}L^2_{t}}
\left\|D_x^{-\frac{1}{4}}f_3\right\|_{L^4_{x}L^\infty_{t}}
\left\|D_x^{-\frac{1}{4}}f_4\right\|_{L^4_{x}L^\infty_{t}}
\left\|J_x^{-1}f_6\right\|_{L^\infty_{x,t}}\\
\lesssim&\; N^{-\frac{7}{2}}\> \|f_1\|_{X_{0,\frac{1}{2}+}}\cdots
\|f_6\|_{X_{0,\frac{1}{2}+}} \endaligned\end{equation*} where we use the fact that$|\xi_5|\sim |\xi_1|$ in this case.

\textbf{Estimate in $A_4$}. Since $|\xi_1+\xi_2|\lesssim |\xi_3|$, we always have
\begin{equation}\label{3.32}
|\bar{M}_6|\lesssim |\xi_1||\xi_2||\xi_3|.
\end{equation}
If $|\xi_j|=|\xi_1|+o(|\xi_1|)$, for all $j=1,\cdots, 5$, then
$|\xi_6|=|\xi_1+\cdots+\xi_5|\gtrsim |\xi_1|$, a contradiction to $A_4$. So we have $|\xi_1|-|\xi_5|\gtrsim |\xi_1|$, which leads to $|\xi_1|\sim |\xi_1-\xi_5|$. In addition, $|\xi_3|\sim |\xi_4|$. Therefore, for any $s>3/8$, as in \eqref{eq-27}, by Lemma \ref{le:strichartz-esimate-3}, we have
\begin{equation*}\aligned
\text{LHS of } \eqref{3.11}\lesssim&\; N^{5s-5}\displaystyle \int dt\int_{A_4}
\frac{\widehat{f_1}(\xi_1)\cdots \widehat{f_6}(\xi_6)}
{|\xi_1|^{s-1} |\xi_2|^{s-1}|\xi_3|^{s-1} |\xi_4|^{s}|\xi_5|^{s}
\langle\xi_6\rangle}\\
\lesssim&\;
N^{-\frac{7}{2}}\left\|D_xf_1f_5\right\|_{L^2_{x,t}}
\left\|D_xf_2\right\|_{L^\infty_{x}L^2_{t}}
\left\|D_x^{-\frac{1}{4}}f_3\right\|_{L^4_{x}L^\infty_{t}}
\left\|D_x^{-\frac{1}{4}}f_4\right\|_{L^4_{x}L^\infty_{t}}
\left\|J_x^{-1}f_6\right\|_{L^\infty_{x,t}}\\
\lesssim&\; N^{-\frac{7}{2}}\> \|f_1\|_{X_{0,\frac{1}{2}+}}\cdots
\|f_6\|_{X_{0,\frac{1}{2}+}}. \endaligned\end{equation*}

\textbf{Estimate in $A_5$}. Again by (\ref{3.32}), Lemma \ref{le:strichartz-esimate-1} and the fact that $|\xi_3|\sim |\xi_4|$, for $s>3/8$, we have
\begin{equation*}\aligned
\text{LHS of }\eqref{3.11}\lesssim&\; N^{6s-6}\displaystyle\int dt\int_{A_5}
\frac{\widehat{f_1}(\xi_1)\cdots \widehat{f_6}(\xi_6)}
{|\xi_1|^{s-1} |\xi_2|^{s-1}|\xi_3|^{s-1} |\xi_4|^{s}|\xi_5|^{s}
|\xi_6|^{s}}\\
\lesssim&\; N^{-4}\left\|D_xf_1\right\|_{L^\infty_{x}L^2_{t}}
\left\|D_xf_2\right\|_{L^\infty_{x}L^2_{t}}
\left\|D_x^{-\frac{1}{4}}f_3\right\|_{L^4_{x}L^\infty_{t}} \cdots
\left\|D_x^{-\frac{1}{4}}f_6\right\|_{L^4_{x}L^\infty_{t}}\\
\lesssim&\; N^{-4}\> \|f_1\|_{X_{0,\frac{1}{2}+}}\cdots
\|f_6\|_{X_{0,\frac{1}{2}+}}. \endaligned\end{equation*} This completes the proof for Proposition \ref{6-linear}.
\end{proof}

Now let us turn to establishing the 10-linear estimate. We first establish a pointwise bound on $\bar{M}_{10}$.
\begin{lem}\label{le:10-linear-bound} If $|\xi_A|\sim |\xi_B|\gtrsim N \gg |\xi_C|\ge |\xi_D|\ge \cdots\ge |\xi_J|$, we have
\begin{equation}
 |\bar{M}_{10}|\lesssim |\xi_C|.
\end{equation}
\end{lem}
\begin{proof}
By Remark \ref{re:two-reductions}, we may assume that $$|\xi_1|\geq \cdots \geq|\xi_{10}|, \text{ and }|\xi_1|\sim |\xi_2|\gtrsim N. $$
Set $\Xi=\big\{3,4,\cdots,10\}$. We rewrite
\begin{equation*}\aligned
\bar{M}_{10} =&\; C\sum\limits_{\{a,\cdots,h\}\subset\Xi}
\Big\{\big[\sigma_6(\xi_1,\xi_2,\xi_a,\xi_b,\xi_c, \xi_d+\xi_e+\xi_f+\xi_g+\xi_h)\\
&\quad +\tilde{\sigma}_6(\xi_1,\xi_2,\xi_a,\xi_b,\xi_c, \xi_d+\xi_e+\xi_f+\xi_g+\xi_h)\big]
(\xi_d+\xi_e+\xi_f+\xi_g+\xi_h)\\
&+\big[\sigma_6(\xi_a,\cdots,\xi_e, \xi_1+\xi_2+\xi_f+\xi_g+\xi_h)\\
&\quad +\tilde{\sigma}_6(\xi_a,\cdots,\xi_e,
\xi_1+\xi_2+\xi_f+\xi_g+\xi_h)\big](\xi_1+\xi_2+\xi_f+\xi_g+\xi_h)\\
&+ \big[\sigma_6(\xi_1,\xi_a,\cdots,\xi_d, \xi_2+\xi_e+\cdots+\xi_h)
(\xi_2+\xi_e+\cdots+\xi_h)\\
&\quad+\sigma_6(\xi_2,\xi_a,\cdots,\xi_d, \xi_1+\xi_e+\cdots+\xi_h)
(\xi_1+\xi_e+\cdots+\xi_h)\big]\\
&+ \big[\tilde{\sigma}_6(\xi_1,\xi_a,\cdots,\xi_d,
\xi_2+\xi_e+\cdots+\xi_h)
(\xi_2+\xi_e+\cdots+\xi_h)\\
&\quad+\tilde{\sigma}_6(\xi_2,\xi_a,\cdots,\xi_d,
\xi_1+\xi_e+\cdots+\xi_h) (\xi_1+\xi_e+\cdots+\xi_h)\big]
\Big\}\\
\equiv&\; \bar{M}_{10}^0+\bar{M}_{10}^1+\bar{M}_{10}^2+\bar{M}_{10}^3.
\endaligned\end{equation*}
Since $|\sigma_6|, |\tilde{\sigma}_6|\lesssim 1$ and $|\xi_d+\cdots+\xi_h|\lesssim |\xi_C|$, it follows that
$$|\bar{M}_{10}^0| \lesssim |\xi_C|.$$
Similarly $|\bar{M}_{10}^1|\lesssim |\xi_3|$, as $|\sigma_6|,|\tilde{\sigma}_6| \lesssim 1$ and
$|\xi_1+\xi_2+\xi_f+\cdots+\xi_h|\lesssim |\xi_3|$.

In order to prove $|\bar{M}_{10}^2|\lesssim |\xi_3|$, by the definition of $\sigma_6$, we only need to show
\begin{equation*}
\begin{split}
&|m(\xi_1)m(\xi_2+\xi_e+\cdots+\xi_h)
(\xi_2+\xi_e+\cdots+\xi_h)\\
&\quad+m(\xi_2)m(\xi_1+\xi_e+\cdots+\xi_h) (\xi_1+\xi_e+\cdots+\xi_h)| \lesssim |\xi_3|.
\end{split}
\end{equation*}
But it follows from the usual mean value theorem twice.

For $\bar{M}_{10}^3$, since $\tilde{\sigma}_6=-\chi_{\Omega}\dfrac{M_6^1}{\alpha_6}-\chi_{\Gamma_6} \sigma_6$, it is concluded by the following lemma and the estimate on $\bar{M}_{10}^2$.
\end{proof}

\begin{lem}
Assume that $|\xi_1|\sim|\xi_2|\gtrsim N \gg |\xi_a|, \cdots|\xi_h|$.
then the following estimates hold:
\begin{align}
\label{eq-28} & \chi_{\Omega_1}(\xi_1,\xi_a,\cdots,\xi_d, \xi_2+\xi_e+\cdots+\xi_h)
           =\chi_{\Omega_1}(\xi_2,\xi_a,\cdots,\xi_d, \xi_1+\xi_e+\cdots+\xi_h);\\
\label{eq-29} &\chi_{\Omega_1}(\xi_1,\xi_a,\cdots,\xi_d, \xi_2+\xi_e+\cdots+\xi_h)
            \cdot\biggl[\dfrac{M_6^1(\xi_1,\xi_a,\cdots,\xi_d, \xi_2+\xi_e+\cdots+\xi_h)}
            {\alpha_6(\xi_1,\xi_a,\cdots,\xi_d, \xi_2+\xi_e+\cdots+\xi_h)} \nonumber \\
           &\times (\xi_2+\xi_e+\cdots+\xi_h)
           +\dfrac{M_6^1(\xi_2,\xi_a,\cdots,\xi_d, \xi_1+\xi_e+\cdots+\xi_h)}
            {\alpha_6(\xi_2,\xi_a,\cdots,\xi_d, \xi_1+\xi_e+\cdots+\xi_h)}
           \cdot (\xi_1+\xi_e+\cdots+\xi_h)\biggl]\\
           &\qquad \lesssim |\xi_3|. \nonumber
\end{align}
\end{lem}
\begin{proof}
By Remark \ref{re:two-reductions}, we may also assume that $$|\xi_1|\geq \cdots \geq|\xi_{10}|, \text{ and }|\xi_1|\sim |\xi_2|\gtrsim N. $$
We first note that $\Omega=\Omega_1$ if $|\xi_1|\sim |\xi_2|\gtrsim N\gg |\xi_3|$. Then we set up some notations for short: for $\{a,b,c,d,e,f,g,h\}\in \Xi$ with $\Xi$ being defined in Lemma \ref{le:10-linear-bound},
\begin{equation*}
\aligned
\overline{\xi_1}&:=\xi_1+\xi_e+\cdots+\xi_h;\\
\overline{\xi_2}&:=\xi_2+\xi_e+\cdots+\xi_h;\\
M&:=M_6^1(\xi_1,\xi_a,\cdots,\xi_d,\xi_2+\xi_e+\cdots+\xi_h);\\
M'&:= M_6^1(\xi_2,\xi_a,\cdots,\xi_d, \xi_1+\xi_e+\cdots+\xi_h);\\
\alpha&:=\alpha_6(\xi_1,\xi_a,\cdots,\xi_d,
\xi_2+\xi_e+\cdots+\xi_h);\\
\alpha'&:=\alpha_6(\xi_2,\xi_a,\cdots,\xi_d,\xi_1+\xi_e+\cdots+\xi_h).
\endaligned
\end{equation*}

To prove \eqref{eq-28}, noting that $|\xi_1|\sim|\xi_2|\sim|\overline{\xi_1}|\sim|\overline{\xi_2}| \gtrsim N \gg |\xi_a|, \cdots|\xi_d|$, we have
$$
|\xi_1^3+\overline{\xi_2}^3|\\
\sim |\xi_1+\overline{\xi_2}|(\xi_1^2+\xi_2^2) \sim |\xi_2^3+\overline{\xi_1}^3|.
$$
Thus \eqref{eq-28} follows from the definition.

To prove \eqref{eq-29}, we need to show that,
\begin{equation}\label{eq-30}
\left|\dfrac{M}{\alpha}\overline{\xi_2}+\dfrac{M'}{\alpha'}\overline{\xi_1}\right|
\lesssim |\xi_3|,\,\forall \,(\xi_1,\xi_2,\xi_a,\cdots,
\xi_h) \in \Omega_1(\xi_1,\xi_a,\cdots,\xi_d, \overline{\xi_2}).
\end{equation}
A computation gives that
\begin{equation}\label{eq-31}
\dfrac{M}{\alpha}\overline{\xi_2}+\dfrac{M'}{\alpha'}\overline{\xi_1}
= \dfrac{\alpha'-\alpha}{\alpha'}\dfrac{M}{\alpha}\overline{\xi_2}
+\dfrac{M-M'}{\alpha'}\overline{\xi_2}
+\dfrac{M'}{\alpha'}(\overline{\xi_1}+\overline{\xi_2}) \equiv
I_1+I_2+I_3.
\end{equation}
For $I_3$, by Lemma \ref{le:6-linear-bound}, we have $|M|\lesssim |\alpha|$ and $|M'|\lesssim |\alpha'|$. Therefore,
\begin{equation}\label{eq-32}
|I_3|\lesssim |\overline{\xi_1}+\overline{\xi_2}| \lesssim |\xi_3|.
\end{equation}
For $I_1$, since $\xi_1-\xi_2-(\overline{\xi_1}-\overline{\xi_2})=0$, the identity \eqref{eq:arithmetic} gives that
\begin{equation*}
\alpha-\alpha'= i\bigl(\xi_1^3+\overline{\xi_2}^3+(-\overline{\xi_1})^3+(-\xi_2)^3\bigr) =
3i(\xi_1+\overline{\xi_2})(\xi_1-\overline{\xi_1})(\xi_1-\xi_2),
\end{equation*} which, together with that fact that $\xi_1\cdot \xi_2 <0$, in turn gives that
 $$|\alpha-\alpha'|\sim |\xi_1+\overline{\xi_2}||\overline{\xi_1}-\xi_1||\xi_1|.$$
Therefore, we see that
\begin{equation}\label{eq-33}
|I_1|\lesssim |\overline{\xi_1}-\xi_1|\lesssim |\xi_3|,
\end{equation}
since $|M|\lesssim|\alpha|$ and the definition of $\Omega_1$ gives that
$$|\alpha'|\gtrsim |\xi_2^3+\overline{\xi_1}^2|\gtrsim |\xi_1+\overline{\xi_2}|(\xi_1^2+\xi_2^2).$$
For $I_2$, by the double mean value theorem in Lemma \ref{le:mvt}, we have
\begin{equation*}
\aligned
|M-M'| =&\; \dfrac{1}{6}|m^2(\xi_1)\xi_1^3+m^2(\overline{\xi_2})\overline{\xi_2}^3
-m^2(\overline{\xi_1})\overline{\xi_1}^3-m^2(\xi_2)\xi_2^3|\\
\lesssim  &\; |\xi_1+\overline{\xi_2}||\overline{\xi_1}-\xi_1||\xi_1|
\endaligned\end{equation*} as $|\xi_1+\overline{\xi_2}| \ll |\xi_1|$ and $|\overline{\xi_1}-\xi_1|\ll |\xi_1|$. Therefore as for $I_1$, we have
$$|I_2| \lesssim |\overline{\xi_1}-\xi_1|\lesssim |\xi_3|. $$
This completes the proof for this lemma.
\end{proof}

Finally the following Proposition will establish the second half of Theorem \ref{thm:main-2}.
\begin{prop}\label{10-linear}
For any $s>1/5$, we have
\begin{equation}
 \left|\displaystyle\int_0^\delta\Lambda_{10}(\bar{M}_{10})\,dt\right|
 \lesssim
 N^{-\frac{15}{4}+}\>
 \|Iu\|_{X_{1,\frac{1}{2}+}^\delta}^{10}.
 \label{Lambda10}
\end{equation}
\end{prop}
\begin{proof}
By Remark \ref{re:two-reductions} or as in the proof of the previous Proposition \ref{6-linear}, we may assume that $$|\xi_1|\ge |\xi_2|\ge |\xi_3| \cdots \geq|\xi_{10}|,\,|\xi_i|\sim N_i, i=1,\cdots,10, \text{ and } N_1\sim N_2\gtrsim N.  $$
By the argument at the beginning of the proof of Proposition \ref{6-linear} and Plancherel's theorem, it suffices to show
\begin{equation}\label{3.13}
\displaystyle \int dt \int_{\Gamma_{10}}\!\!
\frac{\left|\bar{M}_{10}(\xi_1,\cdots,
\xi_{10})\widehat{f_1}(\xi_1) \cdots \widehat{f_{10}}(\xi_{10})\right|} {\langle\xi_1\rangle
m(\xi_1)\cdots \langle\xi_{10}\rangle m(\xi_{10})} \lesssim N^{-\frac{15}{4}+}\>\|f_1\|_{X_{0,\frac{1}{2}+}}\cdots
\|f_{10}\|_{X_{0,\frac{1}{2}+}}.
\end{equation}
We may also assume that the spatial Fourier transforms of all $f_i$ are non-negative. Now we divide it into two regions:
\begin{equation*}\aligned
  B_1 =& \{(\xi_1,\cdots,\xi_{10})\in \Gamma_{10}: |\xi_1|\sim|\xi_2|\gtrsim N\gg |\xi_3|.\}; \\
  B_2 =& \{(\xi_1,\cdots,\xi_{10})\in \Gamma_{10}: |\xi_1|\sim|\xi_2|\geq |\xi_3|\gtrsim N\}.
\endaligned\end{equation*}

\textbf{Estimate in $B_1$}. By Lemma \ref{le:10-linear-bound}, there holds that $|\bar{M}_{10}|\lesssim|\xi_3|$. Then, by Lemmas \ref{le:strichartz-esimate-2} and  \ref{le:strichartz-esimate-3}, by the same reasoning as proving \eqref{eq-27}, the left-hand side of (\ref{3.13}) is bounded by
\begin{equation*}\aligned
N^{2s-2}&  \int dt \int_{B_1} \frac{\widehat{f_1}(\xi_1) \cdots
\widehat{f_{10}}(\xi_{10})} {|\xi_1|^{s} |\xi_2|^{s}
\langle\xi_4\rangle \cdots \langle\xi_{10}\rangle}\\
\lesssim&\; N^{-4}\left\|D_xf_1f_3\right\|_{L^2_{x,t}}
\left\|D_xf_2f_4\right\|_{L^2_{x,t}}
\left\|J_x^{-1}f_5\right\|_{L^\infty_{x,t}} \cdots
\left\|J_x^{-1}f_{10}\right\|_{L^\infty_{x,t}}\\
\lesssim&\; N^{-4}\>\|f_1\|_{X_{0,\frac{1}{2}+}}\cdots
\|f_{10}\|_{X_{0,\frac{1}{2}+}}. \endaligned\end{equation*}

\textbf{Estimate in $B_2$}. By Lemma \ref{le:6-linear-bound}, we have
$|\bar{M}_{10}|\lesssim|\xi_1|$. Then, by the fact
$$
\langle\xi\rangle m(\xi)=\langle\xi\rangle,\text{ for }|\xi|\leq
N;\quad \langle\xi\rangle m(\xi)\sim N^{1-s}|\xi|^s, \text{ for }|\xi|\gtrsim N,
$$
and for $s>1/5$, the left-hand side of (\ref{3.13}) is bounded by
\begin{equation*}\aligned
N^{3s-3}&  \int dt \int_{B_2}
\frac{|\xi_1||\xi_2||\xi_3|^{-\frac{1}{4}}\cdots
|\xi_6|^{-\frac{1}{4}}\> \langle\xi_7\rangle^{-\frac{1}{2}-}\cdots
\langle\xi_{10}\rangle^{-\frac{1}{2}-} \widehat{f_1}(\xi_1)
\cdots \widehat{f_{10}}(\xi_{10})} {|\xi_1|^{2s+1}
|\xi_3|^{s-\frac{1}{4}}\>
m(\xi_4)\langle\xi_4\rangle^{\frac{3}{4}}\cdots
m(\xi_6)\langle\xi_6\rangle^{\frac{3}{4}}\>
m(\xi_7)\langle\xi_7\rangle^{\frac{1}{2}-}\cdots
m(\xi_{10})\langle\xi_{10}\rangle^{\frac{1}{2}-}}\\
\lesssim&
N^{-\frac{15}{4}}\left\|D_xf_1\right\|_{L^\infty_{x}L^2_{t}}
\left\|D_xf_2\right\|_{L^\infty_{x}L^2_{t}}
\left\|D_x^{-\frac{1}{4}}f_3\right\|_{L^4_{x}L^\infty_{t}} \cdots
\left\|D_x^{-\frac{1}{4}}f_6\right\|_{L^4_{x}L^\infty_{t}}\>\\
& \qquad \times
\left\|J_x^{-\frac{1}{2}-}f_7\right\|_{L^\infty_{x,t}} \cdots
\left\|J_x^{-\frac{1}{2}-}f_{10}\right\|_{L^\infty_{x,t}}\\
\lesssim& N^{-\frac{15}{4}}\>\|f_1\|_{X_{0,\frac{1}{2}+}}\cdots
\|f_{10}\|_{X_{0,\frac{1}{2}+}}, \endaligned\end{equation*} where we have used the fact that $|\xi_1|\sim|\xi_2| \geq|\xi_3| \gtrsim N$ and the fact that $s>1/5$.
\end{proof}

\section{Theorem \ref{thm:main-2} implies Theorem \ref{thm:main}}\label{sec:final-argument}
In this section, we show how Theorem \ref{thm:main} is implied by Theorem \ref{thm:main-2} with the help of the modified local theory in Proposition \ref{prop:modified-local}.
\begin{proof}[Proof of Theorem \ref{thm:main}] Fix $u,\,u_0,\, T$ and $3/8<s<1$ as in Theorem \ref{thm:main} and write $A:=1+\|u_0\|_{H^s(\R)}$. Our goal is to show that the solution $u$ to \eqref{gkdv} and \eqref{1.2} exists on $[0,T]$. To this end, we will follow the steps in \cite{CKSTT-03-KDV}. Alternatively, one can also follow the steps in \cite[Section 2]{CKSTT-08} to establish \emph{a priori} bound of the following form
$$ \|u(T)\|_{H^s}\le C(s, \|u_0\|_{H^s}, T).$$
provided that Theorem \ref{thm:main-2} holds true. Then Theorem \ref{thm:main} follows from the local wellposedness theory.

\textbf{Rescaling}. We choose the rescaling parameter $\lambda>1$ which will be determined shortly. We rescale $u$ by
\begin{equation}\label{eq-34}
u_\lambda(x,t)=\lambda^{-1/2}u(x/\lambda,t/\lambda^3);\quad
u_{0,\lambda}(x)=\lambda^{-1/2}u_0(x/\lambda).
\end{equation}
Then $u_\lambda(x,t)$ is still the solution of \eqref{gkdv} with the initial data $u(x,0)= u_{0,\lambda}(x)$; and $u(x,t)$ exists on $[0,T]$ if and only if $u_\lambda(x,t)$ exists on $[0,\lambda^3T]$.

Moreover a computation gives that, for $\lambda>1$,
\begin{equation}\label{eq-45}
\|\partial_xIu_{\lambda}(t)\|_{L^2} \lesssim N^{1-s}/\lambda^s\cdot
\|u(t)\|_{H^s}.
\end{equation}
Hence, if we choose $\lambda\sim N^{\frac{1-s}{s}}$ and specify $t=0$, we have
\begin{equation}\label{eq-42}
\|Iu_{0,\lambda}\|_{H^1}\lesssim_A 1.
\end{equation}
An application of Proposition \ref{prop:modified-local} gives two things: if $\delta$ denotes the lifetime of local solution $u$,
\begin{align}
&\label{eq-43}  \delta\sim_A 1,\\
&\label{eq-44} \|Iu_\lambda(t)\|_{X^\delta_{1,1/2+}}\lesssim \|Iu_{0,\lambda}\|_{H^1_x}.
\end{align}

\textbf{Iteration}. We make two observations. Since $m(\xi)\leq 1$, we first see that
\begin{equation}\label{eq-35}
\|Iu_\lambda(t)\|_{L^2_x}\leq \|u_\lambda(t)\|_{L^2_x}=\|u_{0,\lambda}\|_{L^2_x}=\|u_{0}\|_{L^2_x}<\|Q\|_{L^2}.
\end{equation}
Then by the sharp Gagliardo-Nirenberg inequality (\ref{GN}) and \eqref{eq-45},
\begin{equation}\label{eq-36}
\begin{split}
& \|\partial_xIu_\lambda(t)\|_{L^2_x}^2 \sim E_I^1(u_\lambda(t)),\\
& \|Iu_\lambda(t)\|_{H^1_x}^2\sim \|\partial_xIu_\lambda(t)\|_{L^2_x}^2+\|u_0\|^2_{L^2}\sim E_I^1(u_\lambda(t))+\|u_0\|^2_2.
\end{split}
\end{equation}
By (\ref{E2}) and (\ref{dE2}), the rescaled solution satisfies
\begin{equation}\label{eq-36}
E_I^1(u_{\lambda}(t))=E_I^1(u_{0,\lambda})+\Lambda_6(\tilde{\sigma}_6)(0)
+\displaystyle\int_0^t\left(\Lambda_{6}(\bar{M}_6)
+\Lambda_{10}(\bar{M}_{10})\right)\,ds-\Lambda_6(\tilde{\sigma}_6)(t).
\end{equation}
By Theorem \ref{thm:main-2} and \eqref{eq-44}, we have for any
$t\in [0,\delta]$,
\begin{equation}\label{eq-37}
\aligned
E_I^1(u_{\lambda}(t))\leq &E_I^1(u_{0,\lambda})+C_1N^{0-}\|Iu_{0,\lambda}\|_{H^1}^6\\
&+C_2N^{-\frac{7}{2}+}\left(
\|Iu_{\lambda}\|_{X_{1,\frac{1}{2}+}^\delta}^{6}+
\|Iu_{\lambda}\|_{X_{1,\frac{1}{2}+}^\delta}^{10}\right)+C_1N^{0-}\|Iu_{\lambda}(t)\|_{H^1}^6
\\
\leq &
E_I^1(u_{0,\lambda})+C_1N^{0-}\|Iu_{0,\lambda}\|_{H^1}^6\\
&+C_3N^{-\frac{7}{2}+}\left( \|Iu_{0,\lambda}\|_{H^1}^{6}+
\|Iu_{0,\lambda}\|_{H^1}^{10}\right)+C_1N^{0-}\|Iu_{\lambda}(t)\|_{H^1}^6.
\endaligned
\end{equation} By \eqref{eq-36}, $\|Iu_{\lambda}(t)\|_{H^1_x}^2\sim E_I^1(u_\lambda(t))+\|u_0\|^2_2$, so if choosing a suitably large $N$, a bootstrap argument yields
\begin{equation}\label{eq-38}
 E_I^1(u_{\lambda} (t))
 \leq 2E_I^1(u_{0,\lambda})+2C_1N^{0-}\|Iu_{0,\lambda}\|_{H^1}^6
+2C_3N^{-\frac{7}{2}+}\left(\|Iu_{0,\lambda}\|_{H^1}^{6}+\|Iu_{0,\lambda}\|_{H^1}^{10}\right).
\end{equation}
Thus if assuming that $2E_I^1(u_{0,\lambda})\le C_A$ and choosing a large $N$, we see that
$$\forall \,t\in[0,\delta], \,E_I^1(u_{\lambda}(t))\leq 2C_A.$$
Then we may extend the lifetime of the local solution $u(t)$ to $t\sim_A 2$; in other words, the lifetime of the local solution remain uniformly of size $1$. Repeating this process $M$ times, we obtain
\begin{equation}\label{eq-39}
\aligned E_I^1(u_{\lambda}(t))\leq &
2E_I^1(u_{0,\lambda})+2C_1N^{0-}\|Iu_{0,\lambda}\|_{H^1}^6 \\
&
+2C_3MN^{-\frac{7}{2}+}\left(\|Iu_{0,\lambda}\|_{H^1}^{6}+\|Iu_{0,\lambda}\|_{H^1}^{10}\right).
\endaligned
\end{equation}
Therefore, $E_I^1(u_{\lambda}(t))\leq 2C_A$ provided $M\lesssim N^{7/2-}$, which implies that the solution $u_\lambda$ exists on $[0,N^{7/2-}]$. Hence, $u$ exists on $[0,\lambda^3T]$ with the relation
\begin{equation}\label{eq-40}
N^{\frac{7}{2}-}\gtrsim \lambda^3T\sim N^{\frac{3(1-s)}{s}}T.
\end{equation}
Thus $T=N^{0+}$ as long as $s>6/13$, which implies Theorem \ref{thm:main} by choosing a large $N$.
\end{proof}

\begin{remark}[A polynomial bound]
From the argument above, there exists  a polynomial bound for the solutions on the $H^s(\R)$ norm. Indeed, for $\lambda^s\sim N^{1-s}$, we have for large $t>0$,
\begin{equation}\label{eq-41}
\|u(t)\|_{H^s}\lesssim \|Iu(t)\|_{H^1}\lesssim \lambda^s
\|Iu_\lambda(\lambda^s t)\|_{H^1} \lesssim\lambda^s\sim N^{1-s}
\lesssim t^{\frac{2s(1-s)}{13s-6}+}.
\end{equation}
\end{remark}

\textbf{Acknowledgements.} C. Miao  and G. Xu  were partly supported by the NSF of China (No.10725102, No.10801015). S. Shao was supported by National Science Foundation under agreement No. DMS-0635607 during the early preparations of this work. Y. Wu was partly supported by Beijing International Center for Mathematical Research. Any opinions, findings and
conclusions or recommendations expressed in this paper are those of the authors and do not reflect necessarily the views of the National Science Foundation.

The first, second and fourth authors would like to thank L. G. Farah for informing us a mistake in a previous version. The second author would like to thank Markus Keel for his explanation on Bourgain's high/low trick in ``Fourier truncation method" and the idea of the ``I-method."


\end{document}